\documentclass[runningheads]{llncs}
\usepackage{amsfonts,amsmath,amssymb}
\usepackage[T1]{fontenc}
\usepackage{graphicx}
\usepackage{cleveref}
\usepackage{mathtools}
\usepackage[disable, colorinlistoftodos,bordercolor=orange,backgroundcolor=orange!20,linecolor=orange,textsize=scriptsize]{todonotes}
\usepackage{phdalgo}
\usepackage{enumerate}
\usepackage{caption}
\usepackage{subcaption}
\usepackage{array}
\usepackage{booktabs}
\usepackage{paralist}

\newtheorem{assumption}{Assumption}

\newcommand{\LC}{C}
\newcommand{\actions}{\mathcal{A}}
\newcommand{\timedelays}{\mathcal{T}}

\newcommand{\ie}{i.e.}
\newcommand{\eg}{e.g.}
\newcommand{\R}{\mathbb{R}}
\newcommand{\lex}{\mathrm{lex}}
\newcommand{\leqlex}{\mathrel{\leq_\lex}}

\newcommand{\geqlex}{\mathrel{\geq_\lex}}

\DeclareMathOperator*{\lexmin}{\mathrm{lexmin}}
\DeclareMathOperator{\relint}{\mathrm{relint}}
\DeclarePairedDelimiter{\scalar}{\langle}{\rangle}

\newcommand{\HL}{\mathrm{HL}}

\usepackage{tikz}
\usepackage{tikz-3dplot}
\usetikzlibrary[calc, fit, backgrounds, patterns]

\newcommand{\Pcal}{\mathcal{P}}

\definecolor{colorJunior}{rgb}{0.4,0.6,0.9}   
\definecolor{colorSenior}{rgb}{0.95,0.7,0.4}  
\definecolor{colorBox}{rgb}{0.85,0.4,0.4}     
\definecolor{colorexit}{rgb}{0.5,0.5,0.5}     

\def\p{2.0} 

\tikzset{place/.style={draw,circle,inner sep=2.5pt,semithick}}
\tikzset{transition/.style={rectangle, thick,fill=black, minimum width=2mm,inner ysep=0.5pt}}
\tikzset{vtransition/.style={rectangle, thick, fill=black, minimum height=2mm, inner xsep=0.5pt}}
\tikzset{token/.style={draw,circle,fill=black!80,inner sep=.35pt}}
\tikzset{arrowPetri/.style={->,>=latex,rounded corners=5pt,semithick}}
\tikzset{subgroup/.style={draw,dashed,rounded corners,inner sep=2pt, minimum width=2cm, minimum height=.75cm}}

\begin{document}
\title{Computing the Congestion Phases of Dynamical Systems with Priorities and Application to Emergency Departments}
\titlerunning{Computing Congestion Phases of Dynamical Systems with Priorities}
\author{Xavier~Allamigeon\orcidID{0000-0002-0258-8018} \and
Pascal~Capetillo\orcidID{0009-0001-3305-5417} \and
St\'ephane~Gaubert\orcidID{0000-0002-2777-9988}}
\authorrunning{X. Allamigeon, P. Capetillo and S. Gaubert}
\institute{INRIA and CMAP, \'Ecole polytechnique, Institut polytechnique de Paris,  CNRS\\ France\\
\email{Firstame.Lastname@inria.fr}}
\maketitle              
\begin{abstract}
Medical emergency departments are complex systems in which patients must be treated according to priority rules based on the severity of their condition.  We develop a model of emergency departments using Petri nets with priorities, described by nonmonotone piecewise linear dynamical systems.  The collection of stationary solutions of such systems forms a ``phase diagram'', in which each phase corresponds to a subset of bottleneck resources (like senior doctors, interns, nurses, consultation rooms, etc.). Since the number of phases is generally exponential in the number of resources, developing automated methods is essential to tackle realistic models. We develop a general method to compute congestion diagrams.  A key ingredient is a polynomial time algorithm to test whether a given ``policy'' (configuration of bottleneck tasks) is achievable by a choice of resources. This is done by reduction to a feasibility problem for an unusual class of lexicographic polyhedra. Furthermore, we show that each policy uniquely determines the system's throughput.  We apply our approach to a case study, analyzing a simplified model of an emergency department from Assistance Publique -- H\^opitaux de Paris.

  \keywords{Performance evaluation\and
    Emergency departments \and Piecewise-linear Dynamics \and Petri nets with priorities \and Polyhedral Computation.}
\end{abstract}
\todo{The abstract should briefly summarize the contents of the paper in
150--250 words.}

\section{Introduction}\label{sec:intro}

\subsubsection{Context: Piecewise Linear Models of Timed Discrete Event Systems.}

The ``max-plus approach'' to discrete event systems
allows the analysis of synchronization and concurrency phenomena by means of piecewise linear dynamics. Initially developed to model synchronization phenomena (\ie, the subclass of timed Petri nets called timed-event graphs~\cite{baccelli1992sync,ccggq99,HOW:05,Komenda2017}), it was later extended to {\em monotone} systems~\cite{CGQ95b,Baccelli2004,gaujal2004optimal}, to account for concurrency (resource sharing) via arbitration by preselection rules. More recently, the framework was adapted to include arbitration by priority rules, first to a special case of road traffic in~\cite{Farhi2011}, later to a general setting~\cite{formats2015}, with an application to the staffing of emergency call centers that was subsequently developed in~\cite{boyet2021}.

The general form of the piecewise dynamical systems arising from the max-plus approach involves \emph{counter functions} $z_i \colon \R_{\geq 0}\to \R_{\geq 0}$; each associated with a type of discrete events in the system, such as the firing of a given transition in a timed Petri net. Every $z_i$ is an increasing function of time \(t\), where the value $z_i(t)$ represents the cumulative number of times the event of type $i$ has occurred up to and including time $t$. The dynamical system over the functions $z_i$ then writes as follows:
\begin{align}\label{eq:tds}
\!\!\!  z_i(t) \!=\! \min_{a\in \actions_i} \bigg( r_i^a + \sum_{\tau\in \timedelays} \sum_{j\in [n]} {(P_\tau^a)}_{ij} z_j(t-\tau ) \bigg), \, i\in [n], \; t\geq \max\timedelays \,. \tag{\(D\)} 
\end{align}
Here, the finite set $\timedelays \subset \{0,1,2,\dots\}$ represents the time delays in the system. For every $i \in [n]$,\footnote{Throughout the paper, we use the notation $[k] \coloneqq \{1, \dots, k\}$} $\actions_i$ is a finite set that corresponds to the possible limiting prerequisites for event $i$. The parameters $r_i^a$ and $(P_\tau^a)_{ij}$ take real values. Notably, the parameters $r_i^a$ are used to model the resources of systems (\eg, staffing, equipment, facilities), and play a central role in its behavior.

Performance evaluation can be addressed by seeking {\em stationary solutions} of the dynamics~\eqref{eq:tds}, \ie, solutions of the form $z_i(t)=u_i + \rho_i t$ for all $i \in [n]$, where $u,\rho\in\R^n$. 
Then, the vector $\rho$ represents the system's {\em throughput} with each entry corresponding to the rate of events occuring for each type. When the system is monotone (\eg, when $(P_\tau^a)_{ij}$ are nonnegative) and has a stoichiometric invariant, the system \eqref{eq:tds} can be reformulated as the dynamic programming equation of a semi-Markov decision process.
Then, $\rho$ is uniquely determined by the parameters of~\eqref{eq:tds}, and a pair $(\rho,u)$ can be computed in polynomial time~\cite{CGQ95b,boyet2021}. 

The study of systems governed by priority rules breaks the monotonicity of the dynamics (\ie, $(P_\tau^a)_{ij}$ may take negative values) and poses an additional challenge for computing stationary solutions.
In fact, understanding the conditions of existence and the properties of stationary solutions of nonmonotone systems 
was already stated as Problem \#34 in the list of open problems in control theory collected in 1999, see~\cite{maxplusblondel}.
Progress was made for the explicit dynamics~\eqref{eq:tds} in~\cite{HSCC}, where we identified sufficient conditions, satisfied by a large family of models, that guarantee stationary solutions to exist, independent of the resource allocation (\ie, the values of the parameters $r_i^a$). In contrast, the computation of the stationary solutions has only been done by hand on a few specific and small-size models~\cite{Farhi2011,boyet2021,boyet2022} and has remained unresolved in general. 

\subsubsection{Contributions.} 
We develop an algorithm that computes the stationary solutions of a timed discrete event system governed by a dynamics of the form~\eqref{eq:tds}.
The algorithm iterates over the set of {\em policies}, each specifying a selection $a \in \actions_i$ for every $i \in [n]$ in the minimum in~\eqref{eq:tds}. It relies on two key ingredients. First, we develop a polynomial time procedure that checks whether a policy is strictly feasible, meaning that it is associated to a stationary regime with a nonidentically zero throughput vector (\Cref{th:strict-feasibility}). This boils down to the study of lexicographic polyhedra, \ie, sets defined by finitely many linear inequalities in the lexicographic sense. Lexicographic polyhedra encompass other already studied generalizations of polyhedral models, such as lexicographic linear programming~\cite{Isermann1982lexprog} or polyhedra defined by mixed systems of strict and non-strict linear inequalities~\cite{Bagnara05}. We show that the feasibility problem of such polyhedra can be decided in polynomial time by reduction to successive linear programs. As a second ingredient, we show that, for generic time delays, the throughput associated to a given policy
is uniquely determined by the resource vector, and we provide an explicit formula for it (\Cref{th-unique}) under a condition that $1$ is a semisimple eigenvalue of the matrix associated to this policy. The semisimplicity condition holds for a broad class of systems~\cite{HSCC}.

This work is motivated by a real case study: the dimensioning of a medical emergency department (ED),
carried out as part of an ongoing project with Assistance Publique -- H\^opitaux de Paris. We build a Petri net model of such an ED and apply our algorithm to generate a congestion phase diagram --- a collection of polyhedral regions that, as a function of the resources, identify the congested components of the system, as well as an explicit expression of the throughput. 

The paper is organized as follows. In~\Cref{sec:model}, we present the Petri net model of the ED and the phase diagram returned by our algorithm. \Cref{subsec:principle} deals with a general overview of the algorithm, and the two ingredients, namely checking strict feasibility  of a policy in polynomial time, and determining the throughput vector, are presented in \Cref{subsec:lex-polyhedron,subsec:throughput}. We conclude with further experimental results and implementation details in~\Cref{sec:experiments}. Proofs and extra material can be found in the appendix.

\subsubsection{Related Works.}
As mentioned, we build on a series of works on max-plus and piecewise-linear approaches to discrete event systems~\cite{baccelli1992sync,CGQ95b,gaujal2004optimal,HOW:05}. More recent developments of this approach include invariant space methods~\cite{Katz2007,DiLoreto2010,Perdon2023} and models of resource sharing~\cite{Goltz2022}. We also refer the reader to~\cite{Komenda2017} for a survey.

The idea of computing phase diagrams by enumerating policies was introduced in~\cite{formats2015}
and later applied in~\cite{boyet2021,boyet2022}. It presents two difficulties:
\begin{inparaenum}[1)] 
\item a number of policies which grows exponentially in the size of the system; 
\item efficiently checking whether a given policy is strictly feasible. 
\end{inparaenum}
Previous methods for checking strict feasibility already relied on polyhedral computations and lexicographic inequalities, but they reduced to an exponential number of (usual) linear inequality feasibility problems (see for instance~\cite{boyet2021,boyet2022}). We provide here the first feasibility check that returns in polynomial time. The challenge posed by the exponential number of policies seems irreducible. Indeed, stationary regimes are a generalization of fixed points of ``positive'' tropical polynomial systems. The existence of such a fixed point was shown to be NP-hard in~\cite{boyet2021}. Hence, we do not expect a general polynomial time algorithm for computing stationary regimes.

A classical theorem states that the dynamic programming equations of Markov decision processes admit an {\em invariant half-line}~\cite{Koh80}, another name for a stationary solution $z(t)= u+ \rho t$. Moreover, the term $\rho$ is uniquely determined. \Cref{th-unique} partly extends this result to the case of ``negative probabilities'', showing that the throughput associated to a given policy is unique provided a semisimplicity condition is satisfied. 

The modeling and staffing of emergency organizations, either emergency call centers, or medical emergency departments, has received much attention. These systems have been analyzed using probabilistic networks, particularly queuing networks, see e.g.~\cite{lecuyer,boeufrobert}. The absence of monotonicity induced by the priority rules places these systems outside the scope of ordinary classes of exactly solvable models. In specific cases, ``scaling limits'' of queuing networks (limits of a family of discrete models with a scaling factor tending to infinity) have been obtained, leading to phase diagrams similar to those from our approach~\cite{boeufrobert}.

\todo[inline]{SG: commented the following: (incorporated more or less). There can be exponentially many phases, and it is an open question to develop a general algorithm to compute these phases. To compute phase diagrams.

creates a bottleneck. Consequently, the collection of stationary solutions forms a ``phase diagram'', allowing us to determine the throughput as a function of numbers of resources.  Similar piecewise-linear models have been developed for other priority-based systems, such as road traffic or call centers. In previous studies, phase diagrams were computed manually, leaving an open question: how can automated methods be developed for this task?}

\section{Petri Net Model of an Emergency Department}\label{sec:model}

\begin{figure}[t!]
  \centering
  \def\tkzscl{.25}
  
  \begin{tikzpicture}[scale=\tkzscl,font=\scriptsize]
  
  \node[place, label={[label distance=0.15*\p]above:$\lambda$}] (p_inflow) at (-2*\p, 0) {};
  \node[vtransition, label={[label distance=-0.5*\p]below:$q_{\textrm{in}}$}] (t_inflow) at ($(p_inflow) + (1*\p, 0)$) {};
  \draw[arrowPetri, color=colorexit] (p_inflow) -- (t_inflow);
  \draw[arrowPetri, color=colorexit] (t_inflow) -- ($(t_inflow) + (0.5*\p, 0.25*\p)$) |- ($(p_inflow) + (0, 1*\p)$) -| ($(p_inflow) + (-1*\p, 0)$) -- (p_inflow);
  \node[subgroup] (admin_group) at ($(t_inflow) + (1*\p, -2*\p)$) {Admin Registration};
  \draw[arrowPetri, color=colorexit] (t_inflow) -| (admin_group.north);
  
  \node[subgroup] (triage_group) at ($(admin_group) + (0, -2*\p)$) {Triage};
  \draw[arrowPetri] (admin_group.south) -- (triage_group.north);
  
  \node[place] (p_box_wait) at ($(triage_group) + (0, -2*\p)$) {};
  \node[transition, label={[right]$q_{\textrm{C}}$}] (t_box) at ($(p_box_wait) + (0, -1.5*\p)$) {};
  \draw[arrowPetri] (triage_group.south) -- (p_box_wait);
  \draw[arrowPetri] (p_box_wait) -- (t_box);
  
  \node[place, label={[above left]$N_{\textrm{C}}$}] (p_boxes) at ($(p_box_wait) + (-3*\p, -1*\p)$) {};
  \draw[arrowPetri, color=colorBox] (p_boxes) -- ($(t_box) + (-1*\p, 0.5*\p)$) -- (t_box);
  \node[place] (p_consult_wait) at ($(t_box) + (0, -1*\p)$) {};

  \begin{scope}[local bounding box=junior_consultation]
  \node[transition, label={[right]$q_{\textrm{JC}}$}] (t_junior_start) at ($(p_consult_wait) + (-4*\p, -1*\p)$) {};
  \node[place, label={[left]$\tau_{\textrm{JC}}$}] (p_junior_consult) at ($(t_junior_start) + (0, -1*\p)$) {};
  \node[transition] (t_junior_consult_end) at ($(p_junior_consult) + (0, -1*\p)$) {};
  \node[place] (p_synch_wait) at ($(t_junior_consult_end) + (0, -1*\p)$) {};
  \node[transition, label={[left]$q_{\textrm{JS}}$}] (t_synch) at ($(p_synch_wait) + (0, -1*\p)$) {};
  \node[place, label={[right]$\tau_{\textrm{JS}}$}] (p_synch) at ($(t_synch) + (0, -1*\p)$) {};
  \node[transition] (t_synch_end) at ($(p_synch) + (0, -1*\p)$) {};
  \end{scope}
  \draw[dashed, rounded corners, colorexit] ($(junior_consultation.north west) + (0,1*\p)$) rectangle ($(junior_consultation.south east) + (0,-0.5*\p)$);
  \node at ($(junior_consultation.north) + (0, 1.5*\p)$) {Junior Cons.};

  \node[place, label={[above left]$N_{\textrm{J}}$}] (p_juniors) at ($(p_junior_consult) + (-2*\p, 0)$) {};

  \draw[arrowPetri, ->>] (p_consult_wait) -| (t_junior_start);
  \draw[arrowPetri] (t_junior_start) -- (p_junior_consult);
  \draw[arrowPetri] (p_junior_consult) -- (t_junior_consult_end);
  \draw[arrowPetri] (t_junior_consult_end) -- (p_synch_wait);
  \draw[arrowPetri] (p_synch_wait) -- (t_synch);
  \draw[arrowPetri] (t_synch) -- (p_synch);
  \draw[arrowPetri] (p_synch) -- (t_synch_end);
  \draw[arrowPetri] (t_box) -- (p_consult_wait);
  
  \draw[arrowPetri, colorJunior] (p_juniors) |- ($(t_junior_start) + (-0.5*\p, 0.5*\p)$) -- (t_junior_start);
  \draw[arrowPetri, colorJunior] (t_synch_end) -- ($(t_synch_end) + (-0.5*\p, -0.5*\p)$) -| (p_juniors); 
  
  \begin{scope}[local bounding box=senior_consultation]
  \node[transition, label={[left]$q_{\textrm{SC}}$}] (t_senior_consult) at ($(p_consult_wait) + (4*\p, -1*\p)$) {};
  \node[place, label={[right]$\tau_{\textrm{SC}}$}] (p_senior_consult) at ($(t_senior_consult) + (0, -1*\p)$) {};
  \node[transition] (t_senior_consult_end) at ($(p_senior_consult) + (0, -1*\p)$) {};
  \end{scope}
  \draw[dashed, rounded corners, colorexit] ($(senior_consultation.north west) + (0, 1*\p)$) rectangle ($(senior_consultation.south east) + (0*\p,-0.5*\p)$);
  \node at ($(senior_consultation.north) + (0*\p, 1.5*\p)$) {Senior Cons.};

  \node[place, label={[above right]$N_{\textrm{S}}$}] (p_seniors) at ($(t_senior_consult) + (2*\p, 2.5*\p)$) {};
  
  \draw[arrowPetri, ->] (p_consult_wait) -| (t_senior_consult);
  \draw[arrowPetri] (t_senior_consult) -- (p_senior_consult);
  \draw[arrowPetri] (p_senior_consult) -- (t_senior_consult_end);

  \draw[arrowPetri, colorSenior] (p_seniors) |- ($(t_senior_consult) + (0.5*\p, 0.5*\p)$) -- (t_senior_consult);
  
  \draw[arrowPetri, ->>, colorSenior] (p_seniors) -| ($(t_synch) + (5.5*\p, 1*\p)$) -- ($(t_synch) + (1*\p, 1*\p)$) -- (t_synch);
  
  \node[place] (p_consult_end) at ($(t_box) + (0, -9*\p)$) {};

  \node[vtransition] (t_exit1) at ($(p_consult_end) + (-2*\p, -1*\p)$) {};
  \node[transition] (t_test1) at ($(p_consult_end) + (0, -2*\p)$) {};
  \node[vtransition] (t_care) at ($(p_consult_end) + (2*\p, -1*\p)$) {};

  \begin{scope}[on background layer]
    \node[
      draw,
      fill=white,                  
      pattern=north east lines,    
      pattern color=gray!60,       
      rounded corners,
      inner sep=8pt,
      fit=(p_consult_end) (t_exit1) (t_test1) (t_care)
    ] {};
  \end{scope}

  \draw[arrowPetri] (t_synch_end) -- ($(t_synch_end) + (0, -0.5*\p)$) |- (p_consult_end);
  \draw[arrowPetri] (t_senior_consult_end) |- (p_consult_end);

  \draw[arrowPetri] (p_consult_end) |- (t_exit1) node[near end, above left] {$\pi_{\textrm{exit1}}$};
  \draw[arrowPetri] (p_consult_end) -- (t_test1) node[near end, right] {$\pi_{\textrm{test1}}$};
  \draw[arrowPetri] (p_consult_end) |- (t_care) node[near end, above right] {$\pi_{\textrm{care}}$};

  \node[subgroup] (care_group) at ($(t_care) + (2*\p, -3*\p)$) {Nursing Care};
  \draw[arrowPetri] (t_care) -| (care_group.north);
  
  \node[place] (p_care_end) at ($(care_group.south) + (0, -1.5*\p)$) {};
  \draw[arrowPetri] (care_group.south) -- (p_care_end);
  \node[vtransition] (t_exit2) at ($(p_care_end) + (2*\p, 0)$) {};
  \draw[arrowPetri] (p_care_end) -- (t_exit2) node[midway, above] {$\pi_{\textrm{exit2}}$};
  
  \node[subgroup] (test_group) at ($(p_care_end) + (-4*\p, 0)$) {Diagnostic Tests};
  \draw[arrowPetri] (t_test1) -- (test_group.north);
  
  \draw[arrowPetri] (p_care_end) -- (test_group.east) node[midway, below] {$\pi_{\textrm{test2}}$};
  
  \node[place] (p_fu_wait) at ($(test_group.south) + (0, -1*\p)$) {};
  \node[transition, label={[label distance=0.5*\p]right:$q_{\textrm{EC}}, \tau_{\textrm{EC}}$}] (t_fu) at ($(p_fu_wait) + (0, -1*\p)$) {};
  \draw[arrowPetri] (test_group.south) -- (p_fu_wait);
  \draw[arrowPetri] (p_fu_wait) -- (t_fu);
  
  \draw[arrowPetri, ->>>, colorSenior] (p_seniors) |- ($(p_inflow) + (-7*\p, 2*\p)$) |- ($(p_fu_wait) + (-1*\p, 0)$) -- (t_fu);
  
  \draw[arrowPetri, colorSenior] (t_senior_consult_end) -- ($(t_senior_consult_end) + (1*\p, -1*\p) $) -- ($(t_senior_consult_end) + (3*\p, -1*\p)$) |- (p_seniors);
  \draw[arrowPetri, colorSenior] (t_synch_end) -- ($(t_synch_end) + (0.5*\p, -0.5*\p)$) -| ($(t_senior_consult_end) + (3*\p, -1*\p)$) |- (p_seniors);
  \draw[arrowPetri, colorSenior] (t_fu) -- ($(t_fu) + (0, -\p)$) -| ($(p_fu_wait) + (7*\p, 0)$) |- (p_seniors);
  
  \draw[arrowPetri, colorBox] (t_exit1) -- ($(t_exit1) + (-5*\p, 0)$) |- (p_boxes);
  \draw[arrowPetri, colorBox] (t_test1) |- ($(t_test1) + (-7*\p, -1*\p)$) |- (p_boxes);
  \draw[arrowPetri, colorBox] (care_group.south) -- ($(care_group.south) + (-0.5*\p, -0.5*\p)$) -- ($(care_group.south) + (-11*\p, -0.5*\p)$) |- (p_boxes);
  
  \end{tikzpicture}
  
  \caption{Petri net model of an emergency department. Parameters prefixed by $\tau$ stand for holding times associated with places, those prefixed by $\pi$ for routing proportions, and those prefixed by $N$ for resources (initial marking).}
  \vspace{-.3cm}
  \label{fig:petri_emergency}
\end{figure}
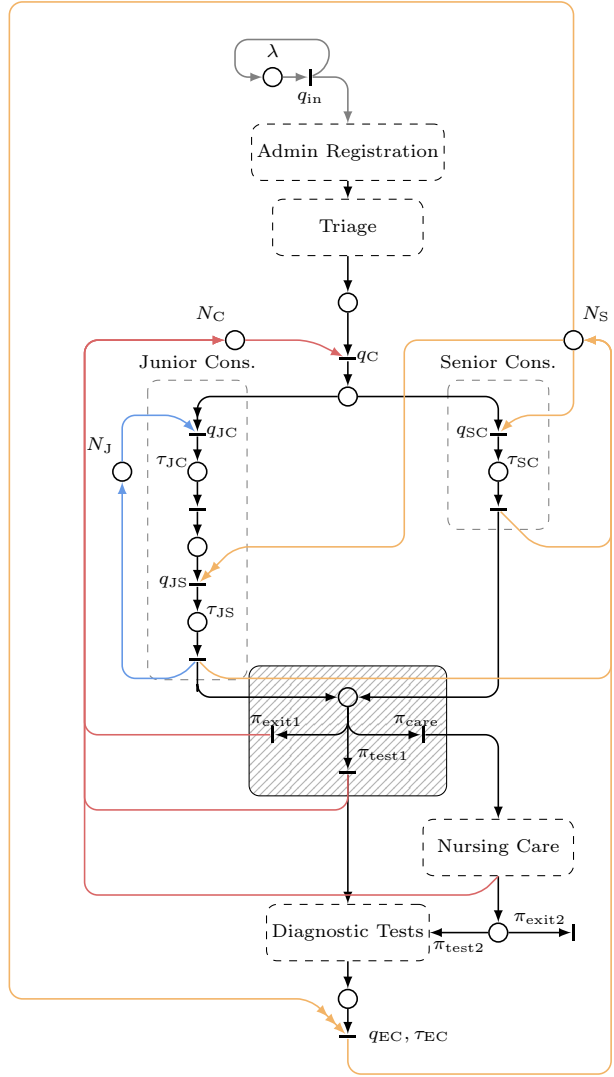

As a motivating example, we develop a Petri net model of an emergency department, illustrated in \Cref{fig:petri_emergency}. This model captures the typical pathways followed by patients within the ED. It consists of the following steps: 
\begin{inparaenum}[(i)] 
\item administrative registration, 
\item triage, 
\item consultation, 
\item nursing care, 
\item diagnostic tests, 
\item final exit consultation.
\end{inparaenum}

The model distinguishes the various types of human resources involved at each stage. Here, administrative registration is handled by medical secretaries, triage by a dedicated pool of nurses or physicians, consultations by junior or senior doctors, nursing care by nurses, and diagnostic tests by specialized technicians. We assume here that the exit consultation is ensured by senior doctors. Additionally, if a patient is initially consulted by a junior doctor, the pathway includes a step in which the junior doctor must validate the diagnosis or treatment with a senior doctor. Finally, certain steps may require material resources; for example, consultations and nursing care are typically conducted in dedicated treatment cubicles. 

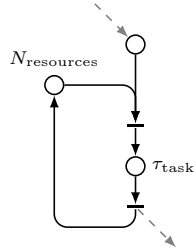
\begin{figure}[h]
  \centering
  \def\tkzscl{.27}
  \vspace{-.2cm}
  
  \begin{tikzpicture}[scale=\tkzscl,font=\scriptsize]
  
  \node[place] (p_wait) at (0, 0) {};
  \node[transition] (t_start) at ($(p_wait) + (0, -2*\p)$) {};
  \node[place] (p_procedure) at ($(t_start) + (0, -1*\p)$) {};
  \node[right of=p_procedure, node distance=.1cm, anchor=west] {$\tau_{\textrm{task}}$};
  \node[transition] (t_end) at ($(p_procedure) + (0, -1*\p)$) {};

  \node[place, label={[above]$N_{\textrm{resources}}$}] (p_resource) at ($(p_wait) + (-2*\p, -1*\p)$) {};

  \draw[arrowPetri] (p_wait) -- (t_start);
  \draw[arrowPetri] (t_start) -- (p_procedure);
  \draw[arrowPetri] (p_procedure) -- (t_end);

  \draw[arrowPetri] (p_resource) -| (t_start);
  \draw[arrowPetri] (t_end) -- ($(t_end) + (0, -0.5*\p)$) -| (p_resource);

  \draw[arrowPetri, dashed, colorexit] ($(p_wait) + (-1*\p, \p)$) -- (p_wait);
  \draw[arrowPetri, dashed, colorexit] (t_end) -- ($(t_end) + (\p,-1*\p)$);

  \end{tikzpicture}

  \caption{Petri net of the administrative registration, triage, nursing, and exam procedures.}
  \vspace{-.3cm}
  \label{fig:petri_subgroup}
\end{figure}

Petri nets provide a natural framework for modeling the synchronization between patients awaiting treatment and the availability of the necessary resources for each step of the care process. Patients, human resources, and material resources are represented by tokens that traverse the Petri net, passing through \emph{places} (depicted by circular nodes) and \emph{transitions} (solid rectangles). Synchronization is modeled by the transitions. To illustrate this, consider the sub-Petri net shown in \Cref{fig:petri_subgroup}, which corresponds to a general task in the ED performed by a single, separate, category of resources. This pattern is replicated in the complete Petri net of \Cref{fig:petri_emergency} to represent administrative registration, triage, nursing care, and diagnostic tests. The place on the left in \Cref{fig:petri_subgroup} represents the pool of available resources for the task, while the place at the top is a waiting area for patients. As soon as each of the two places contains an available token, the top transition fires, consuming two upstream tokens and producing a single token in the downstream place. The latter corresponds to the task itself, and is associated with a holding time $\tau_{\textrm{task}}$. Once this duration has elapsed, the token becomes available to fire the bottom transition, which, in turn, releases the resources back into the pool (via the left arc), and produces another token representing the patient (dashed arc), who proceeds to the next step in the pathway.

Another aspect of emergency departments that Petri nets naturally capture is the concurrent use of human and material resources. In the Petri net, this occurs when multiple transitions compete to consume tokens from a shared upstream place. Our model resolves this in two distinct ways, depending on the modeling needs: \emph{proportion routing} or \emph{priority routing}.

Proportion routing is illustrated by the gray-shaded box at the center of the Petri net in \Cref{fig:petri_emergency}. The shared place collects patients exiting consultations with either junior or senior doctors who then follow one of three possible pathways depending on the severity of their condition: immediate exit, nursing care, or diagnostic tests. In the model, the token is routed to one of the three downstream transitions according to predefined proportions: $\pi_{\textrm{exit1}}$ for immediate exit, $\pi_{\textrm{test1}}$ for diagnostic tests, and $\pi_{\textrm{care}}$ for nursing care. These three proportions sum up to $1$. 

Priority routing is illustrated by the circuit involving senior doctors, given by the arcs highlighted in orange in \Cref{fig:petri_emergency}. Senior doctors are assigned to three distinct tasks, ordered by decreasing priority: exit consultation (transition $q_{\textrm{EC}}$), synchronization with junior doctors for diagnostic confirmation ($q_{\textrm{JS}}$), and initial consultation ($q_{\textrm{SC}}$). The level of priority is visually indicated by the number of arrowheads on the arcs leading to each corresponding transition: more arrowheads denote higher priority. This routing scheme specifies when the three downstream transitions are permitted to fire. For instance, transition $q_{\textrm{JS}}$ can only fire if $q_{\textrm{EC}}$ cannot (\ie, no patient is currently waiting for an exit consultation), regardless of whether patients are waiting for an initial consultation. We point out that priority rules can overlap. For example, patients waiting for initial consultation are assigned in priority to junior doctors (transition $q_{\textrm{JC}}$) over senior doctors ($q_{\textrm{SC}}$). This is consistent with the fact that senior doctors are primarily assigned to tasks requiring higher expertise and should only conduct initial consultations as a last resort, typically when the system becomes overcrowded. 

The resources (human or material) of the emergency department are encoded in the \emph{initial marking} of the Petri net, \ie, the number of tokens initially allocated in each place. When nonzero, this marking is indicated next to the corresponding place. In particular, $N_{\textrm{C}}$, $N_{\textrm{J}}$, and $N_{\textrm{S}}$ represent the number of available cubicles, junior doctors, and senior doctors, respectively. The arrival of patients into the ED is modeled by the top component of the Petri net, which consists of a place with initial marking $\lambda$ and a downstream transition that loops back to it, simulating a constant inflow of patients at rate $\lambda$. 

\tdplotsetmaincoords{63}{20}

\begin{table}[htbp]
  \caption{The congestion phase diagram of the ED model of \Cref{fig:petri_emergency}. We set $\pi_{\textrm{cont}} \coloneqq \pi_{\textrm{exam1}} + \pi_{\textrm{care}}\pi_{\textrm{exam2}}$.}
  \label{tab:congestionED}
  \renewcommand{\arraystretch}{1.2}
  \begin{tabular}{@{}m{0.22\textwidth} m{0.34\textwidth} m{0.42\textwidth}@{}}
  \toprule
  \centering \textbf{Bottleneck type} & 
  \multicolumn{1}{c}{\textbf{Phase region}} & 
  \multicolumn{1}{c}{\textbf{Throughput expressions}} \\
  \midrule
  
  \begin{minipage}[t]{\linewidth}
    \centering
    \textbf{Fluid phase}
  \end{minipage}
  &
  
  \begin{minipage}[c]{\linewidth}
    \vspace{0pt} 
    \centering
    \begin{tikzpicture}[tdplot_main_coords, scale=0.23,
  every node/.style={font=\tiny}]
      \draw[->, thick, dotted, color=gray] (0,0,0) -- (6.0,0,0) node[anchor=north east]{$N_{\textrm{J}}/\lambda$};
      \draw[->, thick, dotted, color=gray] (0,0,0) -- (0,6.0,0) node[anchor=east]{$N_{\textrm{S}}/\lambda$};
      \draw[->, thick, dotted, color=gray] (0,0,0) -- (0,0,6.0) node[anchor=south east]{$N_{\textrm{C}}/\lambda$};
      \draw[dotted, gray] (6.0,0,0) -- (6.0,6.0,0) -- (0,6.0,0);
      \draw[dotted, gray] (0,6.0,0) -- (0,6.0,6.0) -- (0,0,6.0);
      \draw[dotted, gray] (0,0,6.0) -- (6.0,0,6.0) -- (6.0,0,0);
      \draw[dotted, gray] (6.0,6.0,0) -- (6.0,6.0,6.0);
      \draw[dotted, gray] (6.0,6.0,6.0) -- (0,6.0,6.0);
      \draw[dotted, gray] (6.0,0,6.0) -- (6.0,6.0,6.0);
      \fill[fill={rgb,255:red,0; green,255; blue,0}, fill opacity=0.2]
      (3.333, 6.000, 4.667) -- 
      (6.000, 1.307, 4.667) -- 
      (6.000, 6.000, 4.667) -- 
      cycle;
      \fill[fill={rgb,255:red,0; green,255; blue,0}, fill opacity=0.2]
      (3.333, 6.000, 4.667) -- 
      (6.000, 1.307, 4.667) -- 
      (3.333, 1.307, 4.667) -- 
      cycle;
      \fill[fill={rgb,255:red,0; green,255; blue,0}, fill opacity=0.2]
      (3.333, 1.307, 6.000) -- 
      (3.333, 6.000, 4.667) -- 
      (3.333, 1.307, 4.667) -- 
      cycle;
      \fill[fill={rgb,255:red,0; green,255; blue,0}, fill opacity=0.2]
      (3.333, 1.307, 6.000) -- 
      (3.333, 6.000, 4.667) -- 
      (3.333, 6.000, 6.000) -- 
      cycle;
      \fill[fill={rgb,255:red,0; green,255; blue,0}, fill opacity=0.2]
      (3.333, 1.307, 6.000) -- 
      (6.000, 1.307, 4.667) -- 
      (3.333, 1.307, 4.667) -- 
      cycle;
      \fill[fill={rgb,255:red,0; green,255; blue,0}, fill opacity=0.2]
      (3.333, 1.307, 6.000) -- 
      (6.000, 1.307, 6.000) -- 
      (6.000, 1.307, 4.667) -- 
      cycle;
      \fill[fill={rgb,255:red,0; green,255; blue,0}, fill opacity=0.2]
      (6.000, 6.000, 6.000) -- 
      (3.333, 6.000, 4.667) -- 
      (6.000, 6.000, 4.667) -- 
      cycle;
      \fill[fill={rgb,255:red,0; green,255; blue,0}, fill opacity=0.2]
      (6.000, 6.000, 6.000) -- 
      (3.333, 6.000, 4.667) -- 
      (3.333, 6.000, 6.000) -- 
      cycle;
      \fill[fill={rgb,255:red,0; green,255; blue,0}, fill opacity=0.2]
      (6.000, 6.000, 6.000) -- 
      (6.000, 1.307, 4.667) -- 
      (6.000, 6.000, 4.667) -- 
      cycle;
      \fill[fill={rgb,255:red,0; green,255; blue,0}, fill opacity=0.2]
      (6.000, 6.000, 6.000) -- 
      (6.000, 1.307, 6.000) -- 
      (6.000, 1.307, 4.667) -- 
      cycle;
      \fill[fill={rgb,255:red,0; green,255; blue,0}, fill opacity=0.2]
      (6.000, 6.000, 6.000) -- 
      (3.333, 1.307, 6.000) -- 
      (3.333, 6.000, 6.000) -- 
      cycle;
      \fill[fill={rgb,255:red,0; green,255; blue,0}, fill opacity=0.2]
      (6.000, 6.000, 6.000) -- 
      (3.333, 1.307, 6.000) -- 
      (6.000, 1.307, 6.000) -- 
      cycle;
      \draw[black, thick] (6.000,6.000,4.667) -- (6.000,1.307,4.667);
      \draw[black, thick, dashed] (6.000,6.000,4.667) -- (3.333,6.000,4.667);
      \draw[black, thick] (6.000,1.307,4.667) -- (3.333,1.307,4.667);
      \draw[black, thick, dashed] (3.333,1.307,4.667) -- (3.333,6.000,4.667);
      \draw[black, thick] (3.333,1.307,4.667) -- (3.333,1.307,6.000);
      \draw[black, thick, dashed] (3.333,6.000,4.667) -- (3.333,6.000,6.000);
      \draw[black, thick] (3.333,6.000,6.000) -- (3.333,1.307,6.000);
      \draw[black, thick] (6.000,1.307,6.000) -- (3.333,1.307,6.000);
      \draw[black, thick] (6.000,1.307,4.667) -- (6.000,1.307,6.000);
      \draw[black, thick] (6.000,6.000,4.667) -- (6.000,6.000,6.000);
      \draw[black, thick] (3.333,6.000,6.000) -- (6.000,6.000,6.000);
      \draw[black, thick] (6.000,6.000,6.000) -- (6.000,1.307,6.000);
    \end{tikzpicture}
  \end{minipage}
    &
    
  \begin{minipage}[c]{\linewidth}
    \vspace{0pt} 
    \scriptsize
    \[
    \begin{array}{rl}
      \rho_{\textrm{C}} &= \lambda \\[\jot]
      \rho_{\textrm{JC}} &= \lambda \\[\jot]
      \rho_{\textrm{SC}} &= 0 \\[\jot]
      \rho_{\textrm{JS}} &= \lambda
    \end{array}
    \]
  \end{minipage}
  \\
  \midrule
  
  \begin{minipage}[t]{\linewidth}
    \centering
    \textbf{Juniors} limit junior cons.
  \end{minipage}
  &
  
  \begin{minipage}[c]{\linewidth}
    \vspace{0pt} 
    \centering
    \begin{tikzpicture}[tdplot_main_coords, scale=0.23,
  every node/.style={font=\tiny}]
      \draw[->, thick, dotted, color=gray] (0,0,0) -- (6.0,0,0) node[anchor=north east]{$N_{\textrm{J}}/\lambda$};
      \draw[->, thick, dotted, color=gray] (0,0,0) -- (0,6.0,0) node[anchor=east]{$N_{\textrm{S}}/\lambda$};
      \draw[->, thick, dotted, color=gray] (0,0,0) -- (0,0,6.0) node[anchor=south east]{$N_{\textrm{C}}/\lambda$};
      \draw[dotted, gray] (6.0,0,0) -- (6.0,6.0,0) -- (0,6.0,0);
      \draw[dotted, gray] (0,6.0,0) -- (0,6.0,6.0) -- (0,0,6.0);
      \draw[dotted, gray] (0,0,6.0) -- (6.0,0,6.0) -- (6.0,0,0);
      \draw[dotted, gray] (6.0,6.0,0) -- (6.0,6.0,6.0);
      \draw[dotted, gray] (6.0,6.0,6.0) -- (0,6.0,6.0);
      \draw[dotted, gray] (6.0,0,6.0) -- (6.0,6.0,6.0);
      \fill[fill={rgb,255:red,36; green,218; blue,0}, fill opacity=0.2]
      (3.333, 1.307, 6.000) -- 
      (0.000, 2.640, 3.333) -- 
      (3.333, 1.307, 4.667) -- 
      cycle;
      \fill[fill={rgb,255:red,36; green,218; blue,0}, fill opacity=0.2]
      (3.333, 1.307, 6.000) -- 
      (0.000, 2.640, 3.333) -- 
      (0.000, 2.640, 6.000) -- 
      cycle;
      \fill[fill={rgb,255:red,36; green,218; blue,0}, fill opacity=0.2]
      (0.000, 6.000, 6.000) -- 
      (3.333, 1.307, 6.000) -- 
      (0.000, 2.640, 6.000) -- 
      cycle;
      \fill[fill={rgb,255:red,36; green,218; blue,0}, fill opacity=0.2]
      (0.000, 6.000, 6.000) -- 
      (3.333, 1.307, 6.000) -- 
      (3.333, 6.000, 6.000) -- 
      cycle;
      \fill[fill={rgb,255:red,36; green,218; blue,0}, fill opacity=0.2]
      (0.000, 6.000, 6.000) -- 
      (0.000, 2.640, 3.333) -- 
      (0.000, 2.640, 6.000) -- 
      cycle;
      \fill[fill={rgb,255:red,36; green,218; blue,0}, fill opacity=0.2]
      (0.000, 6.000, 6.000) -- 
      (0.000, 2.640, 3.333) -- 
      (0.000, 6.000, 3.333) -- 
      cycle;
      \fill[fill={rgb,255:red,36; green,218; blue,0}, fill opacity=0.2]
      (3.333, 6.000, 4.667) -- 
      (3.333, 1.307, 6.000) -- 
      (3.333, 1.307, 4.667) -- 
      cycle;
      \fill[fill={rgb,255:red,36; green,218; blue,0}, fill opacity=0.2]
      (3.333, 6.000, 4.667) -- 
      (3.333, 1.307, 6.000) -- 
      (3.333, 6.000, 6.000) -- 
      cycle;
      \fill[fill={rgb,255:red,36; green,218; blue,0}, fill opacity=0.2]
      (3.333, 6.000, 4.667) -- 
      (0.000, 6.000, 6.000) -- 
      (3.333, 6.000, 6.000) -- 
      cycle;
      \fill[fill={rgb,255:red,36; green,218; blue,0}, fill opacity=0.2]
      (3.333, 6.000, 4.667) -- 
      (0.000, 6.000, 6.000) -- 
      (0.000, 6.000, 3.333) -- 
      cycle;
      \fill[fill={rgb,255:red,36; green,218; blue,0}, fill opacity=0.2]
      (3.333, 6.000, 4.667) -- 
      (0.000, 2.640, 3.333) -- 
      (3.333, 1.307, 4.667) -- 
      cycle;
      \fill[fill={rgb,255:red,36; green,218; blue,0}, fill opacity=0.2]
      (3.333, 6.000, 4.667) -- 
      (0.000, 2.640, 3.333) -- 
      (0.000, 6.000, 3.333) -- 
      cycle;
      \draw[black, thick] (0.000,2.640,3.333) -- (3.333,1.307,4.667);
      \draw[black, thick] (3.333,1.307,4.667) -- (3.333,1.307,6.000);
      \draw[black, thick] (0.000,2.640,6.000) -- (0.000,2.640,3.333);
      \draw[black, thick] (0.000,2.640,6.000) -- (3.333,1.307,6.000);
      \draw[black, thick] (0.000,2.640,6.000) -- (0.000,6.000,6.000);
      \draw[black, thick] (3.333,6.000,6.000) -- (3.333,1.307,6.000);
      \draw[black, thick] (0.000,6.000,6.000) -- (3.333,6.000,6.000);
      \draw[black, thick, dashed] (0.000,2.640,3.333) -- (0.000,6.000,3.333);
      \draw[black, thick, dashed] (0.000,6.000,3.333) -- (0.000,6.000,6.000);
      \draw[black, thick] (3.333,6.000,4.667) -- (3.333,1.307,4.667);
      \draw[black, thick] (3.333,6.000,6.000) -- (3.333,6.000,4.667);
      \draw[black, thick, dashed] (0.000,6.000,3.333) -- (3.333,6.000,4.667);
    \end{tikzpicture}
  \end{minipage}
  &
  
  \begin{minipage}[c]{\linewidth}
    \vspace{0pt} 
    \centering
    \scriptsize
    \[
    \begin{array}{rl}
      \rho_{\textrm{C}} &= \lambda \\[\jot]
      \rho_{\textrm{JC}} &= \frac{N_{J}}{\tau_{\textrm{JC}} + \tau_{\textrm{JS}}} \\[\jot]
      \rho_{\textrm{SC}} &= \lambda - \frac{N_{J}}{\tau_{\textrm{JC}} + \tau_{\textrm{JS}}} \\[\jot]
      \rho_{\textrm{JS}} &= \frac{N_{J}}{\tau_{\textrm{JC}} + \tau_{\textrm{JS}}}
    \end{array}
    \]
  \end{minipage}
  \\
  \midrule
  
  \begin{minipage}[t]{\linewidth}
    \centering
    \textbf{Cubibles} and \textbf{Juniors} limit senior cons.\  and junior cons., respectively
  \end{minipage}
  &
  
  \begin{minipage}[c]{\linewidth}
    \vspace{0pt} 
    \centering
    \begin{tikzpicture}[tdplot_main_coords, scale=0.23,
  every node/.style={font=\tiny}]
      \draw[->, thick, dotted, color=gray] (0,0,0) -- (6.0,0,0) node[anchor=north east]{$N_{\textrm{J}}/\lambda$};
      \draw[->, thick, dotted, color=gray] (0,0,0) -- (0,6.0,0) node[anchor=east]{$N_{\textrm{S}}/\lambda$};
      \draw[->, thick, dotted, color=gray] (0,0,0) -- (0,0,6.0) node[anchor=south east]{$N_{\textrm{C}}/\lambda$};
      \draw[dotted, gray] (6.0,0,0) -- (6.0,6.0,0) -- (0,6.0,0);
      \draw[dotted, gray] (0,6.0,0) -- (0,6.0,6.0) -- (0,0,6.0);
      \draw[dotted, gray] (0,0,6.0) -- (6.0,0,6.0) -- (6.0,0,0);
      \draw[dotted, gray] (6.0,6.0,0) -- (6.0,6.0,6.0);
      \draw[dotted, gray] (6.0,6.0,6.0) -- (0,6.0,6.0);
      \draw[dotted, gray] (6.0,0,6.0) -- (6.0,6.0,6.0);
      \fill[fill={rgb,255:red,72; green,182; blue,0}, fill opacity=0.2]
      (0.000, 6.000, 3.333) -- 
      (0.000, 6.000, 0.000) -- 
      (3.333, 6.000, 4.667) -- 
      cycle;
      \fill[fill={rgb,255:red,72; green,182; blue,0}, fill opacity=0.2]
      (0.000, 2.640, 3.333) -- 
      (3.333, 1.307, 4.667) -- 
      (0.000, 0.000, 0.000) -- 
      cycle;
      \fill[fill={rgb,255:red,72; green,182; blue,0}, fill opacity=0.2]
      (3.333, 1.307, 4.667) -- 
      (0.000, 6.000, 0.000) -- 
      (3.333, 6.000, 4.667) -- 
      cycle;
      \fill[fill={rgb,255:red,72; green,182; blue,0}, fill opacity=0.2]
      (3.333, 1.307, 4.667) -- 
      (0.000, 6.000, 0.000) -- 
      (0.000, 0.000, 0.000) -- 
      cycle;
      \fill[fill={rgb,255:red,72; green,182; blue,0}, fill opacity=0.2]
      (0.000, 2.640, 3.333) -- 
      (0.000, 6.000, 0.000) -- 
      (0.000, 0.000, 0.000) -- 
      cycle;
      \fill[fill={rgb,255:red,72; green,182; blue,0}, fill opacity=0.2]
      (0.000, 2.640, 3.333) -- 
      (0.000, 6.000, 3.333) -- 
      (0.000, 6.000, 0.000) -- 
      cycle;
      \fill[fill={rgb,255:red,72; green,182; blue,0}, fill opacity=0.2]
      (0.000, 2.640, 3.333) -- 
      (3.333, 1.307, 4.667) -- 
      (3.333, 6.000, 4.667) -- 
      cycle;
      \fill[fill={rgb,255:red,72; green,182; blue,0}, fill opacity=0.2]
      (0.000, 2.640, 3.333) -- 
      (0.000, 6.000, 3.333) -- 
      (3.333, 6.000, 4.667) -- 
      cycle;
      \draw[black, thick, dashed] (0.000,6.000,0.000) -- (0.000,6.000,3.333);
      \draw[black, thick] (0.000,6.000,0.000) -- (3.333,6.000,4.667);
      \draw[black, thick] (0.000,6.000,3.333) -- (3.333,6.000,4.667);
      \draw[black, thick] (3.333,1.307,4.667) -- (0.000,2.640,3.333);
      \draw[black, thick] (0.000,0.000,0.000) -- (3.333,1.307,4.667);
      \draw[black, thick] (0.000,0.000,0.000) -- (0.000,2.640,3.333);
      \draw[black, thick] (3.333,6.000,4.667) -- (3.333,1.307,4.667);
      \draw[black, thick] (0.000,0.000,0.000) -- (0.000,6.000,0.000);
      \draw[black, thick] (0.000,6.000,3.333) -- (0.000,2.640,3.333);
    \end{tikzpicture}
  \end{minipage}
  &
  
  \begin{minipage}[c]{\linewidth}
    \vspace{0pt} 
    \centering
    \scriptsize
    \[
    \begin{array}{rl}
      \rho_{\textrm{C}} & = \frac{(N_{\textrm{C}} - N_{J}) \left(\tau_{\textrm{JC}} + \tau_{\textrm{JS}}\right) + N_{J} \tau_{\textrm{SC}}}{\Delta} \\[\jot]
      \rho_{\textrm{JC}} & = \frac{N_{J}}{\tau_{\textrm{JC}} + \tau_{\textrm{JS}}} \\[\jot]
      \rho_{\textrm{SC}} & = \frac{(N_{\textrm{C}} - N_{J}) \left(\tau_{\textrm{JC}} + \tau_{\textrm{JS}}\right) - N_{J} \pi_{\textrm{care}} \tau_{\textrm{care}}}{\Delta} \\[\jot]
      \rho_{\textrm{JS}} & = \frac{N_{J}}{\tau_{\textrm{JC}} + \tau_{\textrm{JS}}} \, , \enspace  \textrm{where we set} \\[\jot] 
\Delta & \coloneqq \begin{multlined}[t]
\pi_{\textrm{care}} \tau_{\textrm{JC}} \tau_{\textrm{care}} + \pi_{\textrm{care}} \tau_{\textrm{JS}} \tau_{\textrm{care}} \\[-.25cm]
+ \tau_{\textrm{JC}} \tau_{\textrm{SC}} + \tau_{\textrm{SC}} \tau_{\textrm{JS}}
 \end{multlined}
    \end{array}
    \]
  \end{minipage}
  \\
  \midrule
  
  \begin{minipage}[t]{\linewidth}
    \centering
    \textbf{Juniors} and \textbf{Seniors} limit junior cons.\ and senior cons., respecticely.
  \end{minipage}
  &
  
  \begin{minipage}[c]{\linewidth}
    \vspace{0pt} 
    \centering
    \begin{tikzpicture}[tdplot_main_coords, scale=0.23,
    every node/.style={font=\tiny}]
      \draw[->, thick, dotted, color=gray] (0,0,0) -- (6.0,0,0) node[anchor=north east]{$N_{\textrm{J}}/\lambda$};
      \draw[->, thick, dotted, color=gray] (0,0,0) -- (0,6.0,0) node[anchor=east]{$N_{\textrm{S}}/\lambda$};
      \draw[->, thick, dotted, color=gray] (0,0,0) -- (0,0,6.0) node[anchor=south east]{$N_{\textrm{C}}/\lambda$};
      \draw[dotted, gray] (6.0,0,0) -- (6.0,6.0,0) -- (0,6.0,0);
      \draw[dotted, gray] (0,6.0,0) -- (0,6.0,6.0) -- (0,0,6.0);
      \draw[dotted, gray] (0,0,6.0) -- (6.0,0,6.0) -- (6.0,0,0);
      \draw[dotted, gray] (6.0,6.0,0) -- (6.0,6.0,6.0);
      \draw[dotted, gray] (6.0,6.0,6.0) -- (0,6.0,6.0);
      \draw[dotted, gray] (6.0,0,6.0) -- (6.0,6.0,6.0);
      \fill[fill={rgb,255:red,109; green,145; blue,0}, fill opacity=0.2]
      (0.000, 0.000, 6.000) -- 
      (0.000, 2.640, 6.000) -- 
      (3.333, 1.307, 6.000) -- 
      cycle;
      \fill[fill={rgb,255:red,109; green,145; blue,0}, fill opacity=0.2]
      (3.333, 1.307, 4.667) -- 
      (0.000, 2.640, 3.333) -- 
      (0.000, 0.000, 0.000) -- 
      cycle;
      \fill[fill={rgb,255:red,109; green,145; blue,0}, fill opacity=0.2]
      (0.000, 2.640, 3.333) -- 
      (0.000, 0.000, 6.000) -- 
      (0.000, 2.640, 6.000) -- 
      cycle;
      \fill[fill={rgb,255:red,109; green,145; blue,0}, fill opacity=0.2]
      (0.000, 2.640, 3.333) -- 
      (0.000, 0.000, 6.000) -- 
      (0.000, 0.000, 0.000) -- 
      cycle;
      \fill[fill={rgb,255:red,109; green,145; blue,0}, fill opacity=0.2]
      (3.333, 1.307, 4.667) -- 
      (0.000, 0.000, 6.000) -- 
      (0.000, 0.000, 0.000) -- 
      cycle;
      \fill[fill={rgb,255:red,109; green,145; blue,0}, fill opacity=0.2]
      (3.333, 1.307, 4.667) -- 
      (0.000, 0.000, 6.000) -- 
      (3.333, 1.307, 6.000) -- 
      cycle;
      \fill[fill={rgb,255:red,109; green,145; blue,0}, fill opacity=0.2]
      (3.333, 1.307, 4.667) -- 
      (0.000, 2.640, 3.333) -- 
      (0.000, 2.640, 6.000) -- 
      cycle;
      \fill[fill={rgb,255:red,109; green,145; blue,0}, fill opacity=0.2]
      (3.333, 1.307, 4.667) -- 
      (0.000, 2.640, 6.000) -- 
      (3.333, 1.307, 6.000) -- 
      cycle;
      \draw[black, thick] (0.000,0.000,6.000) -- (0.000,2.640,6.000);
      \draw[black, thick] (3.333,1.307,6.000) -- (0.000,2.640,6.000);
      \draw[black, thick] (0.000,0.000,6.000) -- (3.333,1.307,6.000);
      \draw[black, thick, dashed] (0.000,2.640,3.333) -- (3.333,1.307,4.667);
      \draw[black, thick, dashed] (0.000,0.000,0.000) -- (0.000,2.640,3.333);
      \draw[black, thick] (0.000,0.000,0.000) -- (3.333,1.307,4.667);
      \draw[black, thick, dashed] (0.000,2.640,6.000) -- (0.000,2.640,3.333);
      \draw[black, thick] (0.000,0.000,0.000) -- (0.000,0.000,6.000);
      \draw[black, thick] (3.333,1.307,6.000) -- (3.333,1.307,4.667);
    \end{tikzpicture}
  \end{minipage}
  &
  
  \begin{minipage}[c]{\linewidth}
    \vspace{0pt} 
    \centering
    \scriptsize
    \[
    \begin{array}{rl}
      \rho_{\textrm{C}} & = \frac{N_{J} \left(\tau_{\textrm{SC}} - \tau_{\textrm{JS}}\right) + N_{S} \left(\tau_{\textrm{JC}} + \tau_{\textrm{JS}}\right)}{\Delta'} \\[\jot]
      \rho_{\textrm{JC}} & = \frac{N_{J}}{\tau_{\textrm{JC}} + \tau_{\textrm{JS}}} \\[\jot]
      \rho_{\textrm{SC}} & = \frac{N_{S} \left(\tau_{\textrm{JC}} + \tau_{\textrm{JS}}\right) - N_{J} \left(\pi_{\textrm{cont}} \tau_{\textrm{EC}} + \tau_{\textrm{JS}}\right)}{\Delta'} \\[\jot]
      \rho_{\textrm{JS}} & = \frac{N_{J}}{\tau_{\textrm{JC}} + \tau_{\textrm{JS}}} \, , \enspace  \textrm{where we set} \\[\jot] 
\Delta' & \coloneqq \begin{multlined}[t]
\pi_{\textrm{cont}} (\tau_{\textrm{JC}} \tau_{\textrm{EC}} + \tau_{\textrm{EC}} \tau_{\textrm{JS}}) \\[-.25cm]
+ \tau_{\textrm{JC}} \tau_{\textrm{SC}} + \tau_{\textrm{SC}} \tau_{\textrm{JS}}
\end{multlined}
    \end{array}
    \]
  \end{minipage}
  \\
  \midrule

  \begin{minipage}[t]{\linewidth}
    \centering
    \textbf{Seniors} limit the release of cubicles (Juniors limit junior cons.).
  \end{minipage}
  &
  
  \begin{minipage}[c]{\linewidth}
    \vspace{0pt} 
    \centering
    \begin{tikzpicture}[tdplot_main_coords, scale=0.23,
  every node/.style={font=\tiny}]
      \draw[->, thick, dotted, color=gray] (0,0,0) -- (6.0,0,0) node[anchor=north east]{$N_{\textrm{J}}/\lambda$};
      \draw[->, thick, dotted, color=gray] (0,0,0) -- (0,6.0,0) node[anchor=east]{$N_{\textrm{S}}/\lambda$};
      \draw[->, thick, dotted, color=gray] (0,0,0) -- (0,0,6.0) node[anchor=south east]{$N_{\textrm{C}}/\lambda$};
      \draw[dotted, gray] (6.0,0,0) -- (6.0,6.0,0) -- (0,6.0,0);
      \draw[dotted, gray] (0,6.0,0) -- (0,6.0,6.0) -- (0,0,6.0);
      \draw[dotted, gray] (0,0,6.0) -- (6.0,0,6.0) -- (6.0,0,0);
      \draw[dotted, gray] (6.0,6.0,0) -- (6.0,6.0,6.0);
      \draw[dotted, gray] (6.0,6.0,6.0) -- (0,6.0,6.0);
      \draw[dotted, gray] (6.0,0,6.0) -- (6.0,6.0,6.0);
      \fill[fill={rgb,255:red,145; green,109; blue,0}, fill opacity=0.2]
      (0.000, 0.000, 6.000) -- 
      (6.000, 0.000, 6.000) -- 
      (0.000, 0.000, 0.000) -- 
      cycle;
      \fill[fill={rgb,255:red,145; green,109; blue,0}, fill opacity=0.2]
      (3.333, 1.307, 4.667) -- 
      (4.667, 1.307, 6.000) -- 
      (3.333, 1.307, 6.000) -- 
      cycle;
      \fill[fill={rgb,255:red,145; green,109; blue,0}, fill opacity=0.2]
      (4.667, 1.307, 6.000) -- 
      (0.000, 0.000, 6.000) -- 
      (3.333, 1.307, 6.000) -- 
      cycle;
      \fill[fill={rgb,255:red,145; green,109; blue,0}, fill opacity=0.2]
      (4.667, 1.307, 6.000) -- 
      (0.000, 0.000, 6.000) -- 
      (6.000, 0.000, 6.000) -- 
      cycle;
      \fill[fill={rgb,255:red,145; green,109; blue,0}, fill opacity=0.2]
      (3.333, 1.307, 4.667) -- 
      (0.000, 0.000, 6.000) -- 
      (3.333, 1.307, 6.000) -- 
      cycle;
      \fill[fill={rgb,255:red,145; green,109; blue,0}, fill opacity=0.2]
      (3.333, 1.307, 4.667) -- 
      (0.000, 0.000, 6.000) -- 
      (0.000, 0.000, 0.000) -- 
      cycle;
      \fill[fill={rgb,255:red,145; green,109; blue,0}, fill opacity=0.2]
      (3.333, 1.307, 4.667) -- 
      (4.667, 1.307, 6.000) -- 
      (6.000, 0.000, 6.000) -- 
      cycle;
      \fill[fill={rgb,255:red,145; green,109; blue,0}, fill opacity=0.2]
      (3.333, 1.307, 4.667) -- 
      (6.000, 0.000, 6.000) -- 
      (0.000, 0.000, 0.000) -- 
      cycle;
      \draw[black, thick] (6.000,0.000,6.000) -- (0.000,0.000,6.000);
      \draw[black, thick] (6.000,0.000,6.000) -- (0.000,0.000,0.000);
      \draw[black, thick] (0.000,0.000,0.000) -- (0.000,0.000,6.000);
      \draw[black, thick, dashed] (4.667,1.307,6.000) -- (3.333,1.307,4.667);
      \draw[black, thick] (3.333,1.307,6.000) -- (4.667,1.307,6.000);
      \draw[black, thick, dashed] (3.333,1.307,6.000) -- (3.333,1.307,4.667);
      \draw[black, thick] (0.000,0.000,6.000) -- (3.333,1.307,6.000);
      \draw[black, thick] (6.000,0.000,6.000) -- (4.667,1.307,6.000);
      \draw[black, thick, dashed] (0.000,0.000,0.000) -- (3.333,1.307,4.667);
    \end{tikzpicture}
  \end{minipage}
  &
  
  \begin{minipage}[c]{\linewidth}
    \vspace{0pt} 
    \centering
    \scriptsize
    \[
    \begin{array}{rl}
      \rho_{\textrm{C}} & = \frac{N_{S}}{\pi_{\textrm{cont}} \tau_{\textrm{EC}} + \tau_{\textrm{JS}}} \\[\jot]
      \rho_{\textrm{JC}} &= \frac{N_{S}}{\pi_{\textrm{cont}} \tau_{\textrm{EC}} + \tau_{\textrm{JS}}} \\[\jot]
      \rho_{\textrm{SC}} &= 0 \\[\jot]
      \rho_{\textrm{JS}} &= \frac{N_{S}}{\pi_{\textrm{cont}} \tau_{\textrm{EC}} + \tau_{\textrm{JS}}} 
    \end{array}
    \]
  \end{minipage}
  \\
  \midrule

  \begin{minipage}[t]{\linewidth}
    \centering
    \textbf{Cubicles} limit junior cons.
  \end{minipage}
  &
  
  \begin{minipage}[c]{\linewidth}
    \vspace{0pt} 
    \centering
    \begin{tikzpicture}[tdplot_main_coords, scale=0.23,
  every node/.style={font=\tiny}]
      \draw[->, thick, dotted, color=gray] (0,0,0) -- (6.0,0,0) node[anchor=north east]{$N_{\textrm{J}}/\lambda$};
      \draw[->, thick, dotted, color=gray] (0,0,0) -- (0,6.0,0) node[anchor=east]{$N_{\textrm{S}}/\lambda$};
      \draw[->, thick, dotted, color=gray] (0,0,0) -- (0,0,6.0) node[anchor=south east]{$N_{\textrm{C}}/\lambda$};
      \draw[dotted, gray] (6.0,0,0) -- (6.0,6.0,0) -- (0,6.0,0);
      \draw[dotted, gray] (0,6.0,0) -- (0,6.0,6.0) -- (0,0,6.0);
      \draw[dotted, gray] (0,0,6.0) -- (6.0,0,6.0) -- (6.0,0,0);
      \draw[dotted, gray] (6.0,6.0,0) -- (6.0,6.0,6.0);
      \draw[dotted, gray] (6.0,6.0,6.0) -- (0,6.0,6.0);
      \draw[dotted, gray] (6.0,0,6.0) -- (6.0,6.0,6.0);
      \fill[fill={rgb,255:red,182; green,72; blue,0}, fill opacity=0.2]
      (0.000, 0.000, 0.000) -- 
      (6.000, 6.000, 0.000) -- 
      (6.000, 0.000, 0.000) -- 
      cycle;
      \fill[fill={rgb,255:red,182; green,72; blue,0}, fill opacity=0.2]
      (0.000, 0.000, 0.000) -- 
      (6.000, 6.000, 0.000) -- 
      (0.000, 6.000, 0.000) -- 
      cycle;
      \fill[fill={rgb,255:red,182; green,72; blue,0}, fill opacity=0.2]
      (3.333, 6.000, 4.667) -- 
      (0.000, 0.000, 0.000) -- 
      (0.000, 6.000, 0.000) -- 
      cycle;
      \fill[fill={rgb,255:red,182; green,72; blue,0}, fill opacity=0.2]
      (3.333, 6.000, 4.667) -- 
      (0.000, 0.000, 0.000) -- 
      (3.333, 1.307, 4.667) -- 
      cycle;
      \fill[fill={rgb,255:red,182; green,72; blue,0}, fill opacity=0.2]
      (6.000, 1.307, 4.667) -- 
      (0.000, 0.000, 0.000) -- 
      (6.000, 0.000, 0.000) -- 
      cycle;
      \fill[fill={rgb,255:red,182; green,72; blue,0}, fill opacity=0.2]
      (6.000, 1.307, 4.667) -- 
      (0.000, 0.000, 0.000) -- 
      (3.333, 1.307, 4.667) -- 
      cycle;
      \fill[fill={rgb,255:red,182; green,72; blue,0}, fill opacity=0.2]
      (6.000, 6.000, 4.667) -- 
      (6.000, 6.000, 0.000) -- 
      (0.000, 6.000, 0.000) -- 
      cycle;
      \fill[fill={rgb,255:red,182; green,72; blue,0}, fill opacity=0.2]
      (6.000, 6.000, 4.667) -- 
      (3.333, 6.000, 4.667) -- 
      (0.000, 6.000, 0.000) -- 
      cycle;
      \fill[fill={rgb,255:red,182; green,72; blue,0}, fill opacity=0.2]
      (6.000, 6.000, 4.667) -- 
      (6.000, 6.000, 0.000) -- 
      (6.000, 0.000, 0.000) -- 
      cycle;
      \fill[fill={rgb,255:red,182; green,72; blue,0}, fill opacity=0.2]
      (6.000, 6.000, 4.667) -- 
      (6.000, 1.307, 4.667) -- 
      (6.000, 0.000, 0.000) -- 
      cycle;
      \fill[fill={rgb,255:red,182; green,72; blue,0}, fill opacity=0.2]
      (6.000, 6.000, 4.667) -- 
      (3.333, 6.000, 4.667) -- 
      (3.333, 1.307, 4.667) -- 
      cycle;
      \fill[fill={rgb,255:red,182; green,72; blue,0}, fill opacity=0.2]
      (6.000, 6.000, 4.667) -- 
      (6.000, 1.307, 4.667) -- 
      (3.333, 1.307, 4.667) -- 
      cycle;
      \draw[black, thick] (6.000,0.000,0.000) -- (6.000,6.000,0.000);
      \draw[black, thick] (6.000,0.000,0.000) -- (0.000,0.000,0.000);
      \draw[black, thick, dashed] (6.000,6.000,0.000) -- (0.000,6.000,0.000);
      \draw[black, thick] (0.000,6.000,0.000) -- (0.000,0.000,0.000);
      \draw[black, thick] (0.000,6.000,0.000) -- (3.333,6.000,4.667);
      \draw[black, thick] (0.000,0.000,0.000) -- (3.333,1.307,4.667);
      \draw[black, thick] (3.333,1.307,4.667) -- (3.333,6.000,4.667);
      \draw[black, thick] (6.000,0.000,0.000) -- (6.000,1.307,4.667);
      \draw[black, thick] (3.333,1.307,4.667) -- (6.000,1.307,4.667);
      \draw[black, thick] (6.000,6.000,0.000) -- (6.000,6.000,4.667);
      \draw[black, thick] (3.333,6.000,4.667) -- (6.000,6.000,4.667);
      \draw[black, thick] (6.000,6.000,4.667) -- (6.000,1.307,4.667);
    \end{tikzpicture}
  \end{minipage}
  &
  
  \begin{minipage}[c]{\linewidth}
    \vspace{0pt} 
    \centering
    \scriptsize
    \[
    \begin{array}{rl}
      \rho_{\textrm{C}} &= \frac{N_{\textrm{C}}}{\pi_{\textrm{care}} \tau_{\textrm{care}} + \tau_{\textrm{JC}} + \tau_{\textrm{JS}}} \\[\jot]
      \rho_{\textrm{JC}} &= \frac{N_{\textrm{C}}}{\pi_{\textrm{care}} \tau_{\textrm{care}} + \tau_{\textrm{JC}} + \tau_{\textrm{JS}}} \\[\jot]
      \rho_{\textrm{SC}} &= 0 \\[\jot]
      \rho_{\textrm{JS}} &= \frac{N_{\textrm{C}}}{\pi_{\textrm{care}} \tau_{\textrm{care}} + \tau_{\textrm{JC}} + \tau_{\textrm{JS}}} 
    \end{array}
    \]
  \end{minipage}
  \\
  \midrule

  \begin{minipage}[t]{\linewidth}
    \centering
    \textbf{Seniors} limit the release of cubicles (Cubicles limit junior cons.).
  \end{minipage}
  &
  
  \begin{minipage}[c]{\linewidth}
    \vspace{0pt} 
    \centering
    \begin{tikzpicture}[tdplot_main_coords, scale=0.23,
  every node/.style={font=\tiny}]
      \draw[->, thick, dotted, color=gray] (0,0,0) -- (6.0,0,0) node[anchor=north east]{$N_{\textrm{J}}/\lambda$};
      \draw[->, thick, dotted, color=gray] (0,0,0) -- (0,6.0,0) node[anchor=east]{$N_{\textrm{S}}/\lambda$};
      \draw[->, thick, dotted, color=gray] (0,0,0) -- (0,0,6.0) node[anchor=south east]{$N_{\textrm{C}}/\lambda$};
      \draw[dotted, gray] (6.0,0,0) -- (6.0,6.0,0) -- (0,6.0,0);
      \draw[dotted, gray] (0,6.0,0) -- (0,6.0,6.0) -- (0,0,6.0);
      \draw[dotted, gray] (0,0,6.0) -- (6.0,0,6.0) -- (6.0,0,0);
      \draw[dotted, gray] (6.0,6.0,0) -- (6.0,6.0,6.0);
      \draw[dotted, gray] (6.0,6.0,6.0) -- (0,6.0,6.0);
      \draw[dotted, gray] (6.0,0,6.0) -- (6.0,6.0,6.0);
      \fill[fill={rgb,255:red,218; green,36; blue,0}, fill opacity=0.2]
      (6.000, 0.000, 0.000) -- 
      (6.000, 0.000, 6.000) -- 
      (0.000, 0.000, 0.000) -- 
      cycle;
      \fill[fill={rgb,255:red,218; green,36; blue,0}, fill opacity=0.2]
      (6.000, 1.307, 6.000) -- 
      (4.667, 1.307, 6.000) -- 
      (6.000, 0.000, 6.000) -- 
      cycle;
      \fill[fill={rgb,255:red,218; green,36; blue,0}, fill opacity=0.2]
      (3.333, 1.307, 4.667) -- 
      (6.000, 0.000, 0.000) -- 
      (6.000, 1.307, 4.667) -- 
      cycle;
      \fill[fill={rgb,255:red,218; green,36; blue,0}, fill opacity=0.2]
      (3.333, 1.307, 4.667) -- 
      (6.000, 0.000, 0.000) -- 
      (0.000, 0.000, 0.000) -- 
      cycle;
      \fill[fill={rgb,255:red,218; green,36; blue,0}, fill opacity=0.2]
      (4.667, 1.307, 6.000) -- 
      (3.333, 1.307, 4.667) -- 
      (0.000, 0.000, 0.000) -- 
      cycle;
      \fill[fill={rgb,255:red,218; green,36; blue,0}, fill opacity=0.2]
      (4.667, 1.307, 6.000) -- 
      (6.000, 0.000, 6.000) -- 
      (0.000, 0.000, 0.000) -- 
      cycle;
      \fill[fill={rgb,255:red,218; green,36; blue,0}, fill opacity=0.2]
      (6.000, 1.307, 6.000) -- 
      (6.000, 0.000, 0.000) -- 
      (6.000, 1.307, 4.667) -- 
      cycle;
      \fill[fill={rgb,255:red,218; green,36; blue,0}, fill opacity=0.2]
      (6.000, 1.307, 6.000) -- 
      (6.000, 0.000, 0.000) -- 
      (6.000, 0.000, 6.000) -- 
      cycle;
      \fill[fill={rgb,255:red,218; green,36; blue,0}, fill opacity=0.2]
      (6.000, 1.307, 6.000) -- 
      (3.333, 1.307, 4.667) -- 
      (6.000, 1.307, 4.667) -- 
      cycle;
      \fill[fill={rgb,255:red,218; green,36; blue,0}, fill opacity=0.2]
      (6.000, 1.307, 6.000) -- 
      (4.667, 1.307, 6.000) -- 
      (3.333, 1.307, 4.667) -- 
      cycle;
      \draw[black, thick] (6.000,0.000,6.000) -- (6.000,0.000,0.000);
      \draw[black, thick] (6.000,0.000,6.000) -- (0.000,0.000,0.000);
      \draw[black, thick] (0.000,0.000,0.000) -- (6.000,0.000,0.000);
      \draw[black, thick] (4.667,1.307,6.000) -- (6.000,1.307,6.000);
      \draw[black, thick] (6.000,0.000,6.000) -- (4.667,1.307,6.000);
      \draw[black, thick] (6.000,0.000,6.000) -- (6.000,1.307,6.000);
      \draw[black, thick] (6.000,0.000,0.000) -- (6.000,1.307,4.667);
      \draw[black, thick, dashed] (6.000,1.307,4.667) -- (3.333,1.307,4.667);
      \draw[black, thick] (0.000,0.000,0.000) -- (3.333,1.307,4.667);
      \draw[black, thick] (3.333,1.307,4.667) -- (4.667,1.307,6.000);
      \draw[black, thick] (6.000,1.307,4.667) -- (6.000,1.307,6.000);
    \end{tikzpicture}
  \end{minipage}
  &
  
  \begin{minipage}[c]{\linewidth}
    \vspace{0pt} 
    \centering
    \scriptsize
    \[
    \begin{array}{rl}
      \rho_{\textrm{C}} &= \frac{N_{S}}{\pi_{\textrm{cont}} \tau_{\textrm{EC}} + \tau_{\textrm{JS}}} \\[\jot]
      \rho_{\textrm{JC}} &= \frac{N_{S}}{\pi_{\textrm{cont}} \tau_{\textrm{EC}} + \tau_{\textrm{JS}}} \\[\jot]
      \rho_{\textrm{SC}} &= 0 \\[\jot]
      \rho_{\textrm{JS}} &= \frac{N_{S}}{\pi_{\textrm{cont}} \tau_{\textrm{EC}} + \tau_{\textrm{JS}}} 
    \end{array}
    \]
  \end{minipage}
  \\
  \midrule
  \bottomrule
  \end{tabular}
\end{table}

\subsubsection*{Congestion Phase Diagram of the ED.} 

Our goal is to assess the performance and congestion levels of an emergency department (ED) based on the available resources. As outlined in \Cref{sec:intro}, we use the throughputs of transitions as a performance metric. Here, they correspond to the rates at which activities, such as consultations or diagnostic tests, are completed. We denote the thoughput vector by \(\rho\), where \(\rho_i\) represents the firing rate of transition \(q_i\).

The congestion phase diagram provides a compact representation of how the performance of the ED evolves under varying resource allocations. It partitions the space of resources into polyhedral cells, each corresponding to a subset of resources that act as bottlenecks, limiting the flow of patients through the system. On each such cell, the throughput vector \(\rho\) is given by an explicit formula determined by the active bottlenecks.

As an illustration, \Cref{tab:congestionED} presents the congestion phase diagram generated by the algorithm developed in the following sections, applied to the Petri net model of the ED. Due to space constraints, the table restricts to the situation where only junior doctors, senior doctors, and cubicles are treated as potentially limiting resources (all other resources are assumed to be infinite, making the corresponding parts of the model fluid). Despite this simplification, these three resources capture the most intricate interactions in the system, particularly those arising from overlapping priority rules. 

The first column of the table identifies the limiting resources and where they act as bottlenecks, the second exhibits the corresponding cell in the resource space, and the third provides the explicit expressions for the throughputs of the relevant transitions. Since the system dynamics are homogeneous (\ie, invariant under a simultaneous scaling of the arrival rate $\lambda$ and the resources by the same factor), we represent every cell using the projective coordinates $(N_{\textrm{J}} / \lambda, N_{\textrm{S}} / \lambda, N_{\textrm{C}} / \lambda)$. 

We refer to \Cref{sec:experiments} for the analysis of the diagram as well as a short discussion of the ED model with all types of resources.

\section{Principle of the Algorithm}\label{sec:max-plus}\label{subsec:principle}

\subsection{Piecewise Linear Dynamics of Petri Nets}

We briefly recall how the piecewise linear dynamics of the form~\eqref{eq:tds} can be  explicitly constructed from a Petri net with priority rules. The counter functions are associated to transitions, \ie, $z_i(t)$ counts the number of times the transition $q_i$ has fired up to and including time $t$. We denote by $n$ the number of transitions in the Petri net. The set $\actions_i$ then corresponds to the places located upstream the transition $q_i$. The parameter $r_i^a$ is then set to the initial marking of the place $a$, \ie, it corresponds to the amount of resources allocated to this place. Holding times in the Petri net are gathered in the set $\timedelays$. The way the parameters  ${(P_\tau^a)}_{ij}$ are specified depends on the kind of routing used from the upstream place $a$. In the case of a proportion routing, ${(P_\tau^a)}_{ij}$ derives from the corresponding proportion. If the upstream place is subject to a priority rule, this parameter can also take $\pm 1$ values. For a more detailed description, we refer to~\cite{formats2015,boyet2021}, where it is proved that an integer version of the dynamical system~\eqref{eq:tds} (in which the counter functions are integer values and routing probabilities are replaced by preselection functions) is in correspondence with the semantics of the Petri net. The continuous setting we consider here, where counter functions collect fractions of tokens and can take real values, has been proved to be asymptotically tight under a scaling limit, for monotone systems~\cite[Chapter~2]{boyet2022}.

As an illustration, we provide in~\Cref{tab:counters} a dynamical system that describes the part of the model of ED restricted to the interactions of patients with junior and senior doctors as well as cubicles (\ie, the dashed components in~\Cref{fig:petri_emergency} are omitted). In this simplified model, we can reduce the dynamics to the five counter functions \(z_{\textrm{C}}\), $z_{\textrm{JC}}$, $z_{\textrm{SC}}$, $z_{\textrm{JS}}$ and $z_{\textrm{EC}}$, corresponding to the transitions $q_{\textrm{C}}$, $q_{\textrm{JC}}$, $q_{\textrm{SC}}$, $q_{\textrm{JS}}$ and $q_{\textrm{EC}}$ respectively. As shown in~\Cref{tab:counters}, the resources of the ED, namely $N_{\textrm{J}}$, $N_{\textrm{S}}$ and $N_{\textrm{C}}$, appears in the ``affine part'' of the dynamics, \ie, in the parameters $r_i^a$, while the duration of tasks, the proportion of patients going to nursing care, exams, etc, are collected in the ``linear part''. The priority rules introduce subtractions of counter functions in some of the linear terms. As an example, the first term in the minimum specifying $z_{\textrm{SC}}(t)$ (third equation) represents the cumulative number of senior doctors that have been available from the senior doctor pool up to time $t$, \ie, $N_{\textrm{S}} + z_{\textrm{JS}} (t - \tau_{\textrm{JS}}) + z_{\textrm{SC}} (t - \tau_{\textrm{SC}}) + z_{\textrm{EC}} (t - \tau_{\textrm{EC}})$, to which we have subtracted the number $z_{\textrm{JS}}(t) + z_{\textrm{EC}}(t)$ of senior doctors that have been already assigned to the higher priority tasks, \ie, exit consultation and diagnosis validation with junior doctors.

\begin{figure}[t]
{\scriptsize
\[ 
\begin{aligned}
  z_{\textrm{C}}(t) &= \min\Big( \lambda t,\ 
    N_{\textrm{C}} + (1 - \pi_{\textrm{care}}) (z_{\textrm{JS}} (t - \tau_{\textrm{JS}}) 
    + z_{\textrm{SC}} (t - \tau_{\textrm{SC}})) \\[-.1cm]
    &\quad + \pi_{\textrm{care}} (z_{\textrm{JS}} (t - \tau_{\textrm{JS}} - \tau_{\textrm{care}})
    + z_{\textrm{SC}} (t - \tau_{\textrm{SC}} - \tau_{\textrm{care}})) \Big) \\[1ex]
  z_{\textrm{JC}}(t) &= \min\Big(
    N_{\textrm{J}} + z_{\textrm{JS}} (t - \tau_{\textrm{JS}}),\ 
    z_{\textrm{C}} (t) - z_{\textrm{SC}} (t^-)
  \Big) \\[1ex]
  z_{\textrm{SC}}(t) &= \min\Big(
    N_{\textrm{S}} + z_{\textrm{JS}} (t - \tau_{\textrm{JS}}) + z_{\textrm{SC}} (t - \tau_{\textrm{SC}}) + z_{\textrm{EC}} (t - \tau_{\textrm{EC}}) \\[-.1cm]
    &\quad - z_{\textrm{JS}} (t) - z_{\textrm{EC}} (t),\ 
    z_{\textrm{C}} (t) - z_{\textrm{JC}} (t)
  \Big) \\[1ex]
  z_{\textrm{JS}}(t) &= \min\Big(
    N_{\textrm{S}} + z_{\textrm{JS}} (t - \tau_{\textrm{JS}}) + z_{\textrm{SC}} (t - \tau_{\textrm{SC}}) + z_{\textrm{EC}} (t - \tau_{\textrm{EC}}) \\[-.1cm]
    &\quad - z_{\textrm{SC}} (t^-) - z_{\textrm{EC}} (t),\ 
    z_{\textrm{JC}} (t - \tau_{\textrm{JC}})
  \Big) \\[1ex]
  z_{\textrm{EC}}(t) &= \min\Big(
    N_{\textrm{S}} + z_{\textrm{JS}} (t - \tau_{\textrm{JS}}) + z_{\textrm{SC}} (t - \tau_{\textrm{SC}}) + z_{\textrm{EC}} (t - \tau_{\textrm{EC}}) \\[-.1cm]
    &\quad - z_{\textrm{SC}} (t^-) - z_{\textrm{JS}} (t^-),\ 
    \pi_{E1} (z_{\textrm{JS}} (t - \tau_{\textrm{JS}} - \tau_{\textrm{test}}) + z_{\textrm{SC}} (t - \tau_{\textrm{SC}} - \tau_{\textrm{test}})) \\[-.1cm]
    &\quad + \pi_{E2} \pi_{\textrm{care}} (z_{\textrm{JS}} (t - \tau_{\textrm{JS}} - \tau_{\textrm{care}} - \tau_{\textrm{test}})
    + z_{\textrm{SC}} (t - \tau_{\textrm{SC}} - \tau_{\textrm{care}} - \tau_{\textrm{test}}))
  \Big)
\end{aligned}
\]}
\caption{The dynamical system over the counter functions of the ED Petri net.}\label{tab:counters}
\end{figure}

\subsection{Stationary Regimes and Lexicographic Constraints}\label{subsec:lex-constraints}

We seek affine stationary regimes of the dynamics~\eqref{eq:tds}, \ie, solutions of form $z_i \colon t \mapsto u_i + \rho_i t$ for all $i \in [n]$. Substituting such soltions in~\eqref{eq:tds} yields a system of lexicographic constraints over the pairs $(\rho_i, u_i)$. Indeed, the terms achieving the minimum in~\eqref{eq:tds} for all $t$ large enough are fully determined by the lexicographic ordering: two affine functions $f(t) = w + \eta t$ and $g(t) = w' + \eta' t$ satisfy $f(t) \leq g(t)$ for all $t$ large enough if and only if $(\eta, w) \leqlex (\eta', w')$, \ie, either $\eta < \eta'$, or $\eta = \eta'$ and $w \leq w'$. In this way, finding an affine stationary solution to~\eqref{eq:tds} can be shown to be equivalent to solving the following system over the variables $(\rho, u) \in \R^n \times \R^n$
\begin{equation}\label{eq:stationary-lex}
\forall i \in [n] \, , \quad (\rho_i, u_i) = \lexmin_{a \in \actions_i} \Big([P^a]_i \rho, r_i^a + [P^a]_i u - [\bar P^a]_i \rho \Big) \, \tag{$\HL$}
\end{equation}
where $\lexmin$ stands for the minimum operator w.r.t.~the lexicographic order~$\leqlex$, $P^a \coloneqq \sum_{\tau \in \timedelays} P^a_\tau$, $\bar P^a \coloneqq \sum_{\tau \in \timedelays} \tau  P^a_\tau$, and $[M]_i$ stands for the $i$th row of the matrix $M$. We refer to~\cite[Lemma~7]{HSCC} for a proof.

As briefly explained in \Cref{sec:intro}, our approach to solve the system~\eqref{eq:tds} relies on policies. Formally, a {\em policy\/} \(\sigma \colon [n]\to \cup_{i\in [n]} \actions_i\) is a map such that \(\sigma(i)\in \actions_i\) for all \(i \in [n]\). If $(\rho, u) \in \R^n \times \R^n$ is a solution of~\eqref{eq:stationary-lex}, then there exists a policy~$\sigma$ such that, for all $i \in  [n]$, 
\begin{equation}\label{eq:policy-lex}
\begin{aligned}
(\rho_i, u_i) & = \Big([P^{\sigma(i)}]_i \rho, \; r_i^{\sigma(i)} + {[P^{\sigma(i)}]}_i u - [\bar P^{\sigma(i)}]_i \rho \Big) \, , \\
(\rho_i, u_i) & \leqlex 
\Big([P^a]_i \rho, \; r_i^a  + {[P^a]}_i u - {[\bar P^a]}_i \rho \Big) 
\quad \text{for all} \enspace a \in \actions_i \, , a \neq \sigma(i) \, .
\end{aligned} \tag{$\HL_\sigma$}
\end{equation}
In this case, the policy determines the nature of the congestion in the Petri net as follows: for each $i \in [n]$, the transition $q_i$ is constrained by the flow through the upstream place associated with $a = \sigma(i) \in \actions_i$. Conversely, any solution to a system of the form~\eqref{eq:policy-lex} is also a solution to~\eqref{eq:stationary-lex}, meaning that it corresponds to an affine stationary regime.

As a consequence, determining the congestion phase diagram amounts to identifying, for every policy $\sigma$, the set $C_\sigma$ of resources $(r_i^a)$ for which the system~\eqref{eq:policy-lex} admits a solution $(\rho, u)$ with $\rho \geq 0$. However, the trivial solution $(\rho_i, u_i) = (0, 0)$ for all $i \in [n]$, corresponding to a fully congested regime, is always feasible regardless of the policy (\eg, setting $r_i^{\sigma(i)} = 0$ for all $i$). We exclude this special case from the analysis and focus instead on \emph{strictly feasible policies}. These are the policies $\sigma$ for which there exists a resource allocation $(r_i^a)$ such that system~\eqref{eq:policy-lex} has a solution $(\rho, u)$ with a throughput vector $\rho$ that is nonnegative and nonzero. The \emph{congestion phase diagram} is then defined as the collection of cells ${C_\sigma}$ where $\sigma$ ranges over the strictly feasible policies. To compute this diagram, we propose to iterate over all policies $\sigma$, and for each one, determine whether it is strictly feasible and, if so, the corresponding cell $C_\sigma$ and the set of stationary regimes defined by the solutions of system~\eqref{eq:policy-lex}. This is the purpose of the techniques developed in \Cref{subsec:lex-polyhedron,subsec:throughput}.

\section{Solving Lexicographic Inequalities in Polynomial Time}\label{subsec:lex-polyhedron}

The systems~\eqref{eq:policy-lex} gives rise to the study of \emph{lexicographic polyhedra} (or \emph{lex-polyhedra} for short), which are defined as the set of solutions $x \in \R^N$ of finitely many \emph{lexicographic linear inequalities} (\emph{lex-inequalities} for short)
\begin{equation}\label{eq:lex-ineq}
\big(\scalar{\alpha^{1}_i, x}, \dots, \scalar{\alpha^{d_i}_i, x}\big) \leqlex \big(b^{1}_i, \dots, b^{d_i}_i\big) \enspace , \quad i \in [m]
\end{equation}
where $\leqlex$ stands for the lexicographic order over $\R^{d_i}$ ($d_i \geq 1$), the $\alpha^{j}_i$ are $N$-vectors ($i \in [m]$, $j \leq d_i$), and $\scalar{y, x} \coloneqq \sum_k y_k x_k$. The integer $d_i$ is called the \emph{depth} of the inequality. Note that depth-$1$ lex-inequalities precisely correspond to usual linear inequalities. Similarly, a usual linear equality can be simply encoded as a pair of opposed depth-$1$ lex-inequalities.

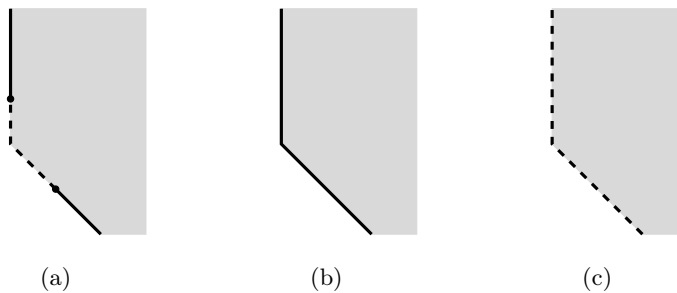
\begin{figure}[t]
\begin{center}
\begin{tikzpicture}[scale=.6]
\begin{scope}
\filldraw[gray!30] (0,3) -- (0,0) -- (2,-2) -- (3,-2) -- (3,3) -- cycle;
\filldraw (0,1) circle (2pt);
\filldraw (1,-1) circle (2pt);
\draw[very thick] (0,1) -- (0,3);
\draw[very thick,dashed] (1,-1) -- (0,0) -- (0,1);
\draw[very thick] (1,-1) -- (2,-2);
\node at (1,-3) {(a)};
\end{scope}

\begin{scope}[xshift=6cm]
\filldraw[gray!30] (0,3) -- (0,0) -- (2,-2) -- (3,-2) -- (3,3) -- cycle;
\draw[very thick] (2,-2) -- (0,0) -- (0,3);
\node at (1,-3) {(b)};
\end{scope}

\begin{scope}[xshift=12cm]
\filldraw[gray!30] (0,3) -- (0,0) -- (2,-2) -- (3,-2) -- (3,3) -- cycle;
\draw[very thick,dashed] (2,-2) -- (0,0) -- (0,3);
\node at (1,-3) {(c)};
\end{scope}
\end{tikzpicture}
\end{center}
\caption{(a) The lex-polyhedron defined by $(x_1, x_2) \geqlex (0, 1)$, $(x_1, x_1) \geqlex (-x_2, 1)$. In the boundary, solid lines and points are included, while dash lines are excluded. (b) The max-front polyhedron $\Pcal^{f^\star}$ and (c) its relative interior. }\label{fig:lex-polyhedron}
\end{figure}

To the best of our knowledge, this way of extending polyhedra to a lexicographic setting has not been studied in the literature. Several works address lexicographic linear programming~\cite{Isermann1982lexprog}, in which only the objective function to be optimized is lexicographic, while the constraint set remains a standard polyhedron. Our lexicographic polyhedra can be seen as a generalization of this framework (sublevel sets of the objective function can be equivalently modeled as one lexicographic inequality), but the fact that multiple inequalities are lexicographic fundamentally changes the nature of the problem. Indeed, although lexicographic polyhedra are still convex sets, they can exhibit particular structures, even in the case of only two lexicographic inequalities, as illustrated in \Cref{fig:lex-polyhedron}(a). Lex-polyhedra also encompass \emph{non-necessarily closed polyhedra}~\cite{Bagnara05} that have been widely applied in abstract interpretation to infer linear inequality program invariants, and that correspond to mixed systems of strict and non-strict linear inequalities. Indeed, a strict linear inequality $\scalar{\alpha, x} < b$ can be equivalently encoded by a depth-$2$ lex-inequality $(\scalar{\alpha,x}, 1) \leqlex (b, 0)$.  

To make the exposition of the algorithm easier, we consider a lex-polyhedron $\Pcal_\lex$ defined by a system of lex-inequalities in \emph{standard form}, \ie,  
\begin{equation}\label{eq:lex-polyhedron}
\begin{aligned}
A x + s & = b \, , \\
\forall i \in [m] \, , (s^1_i, \dots, s^{d_i}_i) & \geqlex 0 \, 
\end{aligned}\tag{$\textsc{Lex-SF}$} 
\end{equation}
where the variables are now $x = (x_1, \dots, x_N)$ and $s = (s_i^j)_{i \in [m] , \, j \in [d_i]}$, $A \in \R^{D \times N}$ and $b \in \R^D$, and $D \coloneqq \sum_{i = 1}^m d_i$ is the total depth. In consequence, $\Pcal_\lex$ is a subset of $\R^{N + D}$. We remark that reformulating systems of lex-inequalities of the form~\eqref{eq:lex-ineq} to the system~\eqref{eq:lex-polyhedron} is simply achieved by introducing the \emph{slack variables} $s_i^j$ in every component of the lex-inequalities, and constraining these groups of variables to be nonnegative in the lexicographic sense. We refer to \Cref{ex:front} at the end of the section for an illustration.

A \emph{front} is a $m$-vector $f$ such that $f_i \in [d_i+1]$ for all $i \in [m]$. We denote $s^{< f}$ the vector consisting of the variables $s^j_i$ for $i \in [m]$ and $j < f_i$, and $s^f$ that of the variables $s^{f_i}_i$ for the $i \in [m]$ such that $f_i \leq d_i$. We set
\[
\Pcal^f \coloneqq \big\{ (x,s) \in \R^{N + D} \colon A x + s = b \, , \; s^{< f} = 0 \, , \; s^f \geq 0 \big\} \, .
\]
The sets $\Pcal^f$ are (closed) convex polyhedra. Our aim is to approximate the lex-polyhedron $\Pcal^\lex$ by means of these polyhedra. Fronts can be ordered componentwise, \ie, $f \leq f'$ if $f_i \leq f'_i$ for all $i \in [m]$. Remark that $\Pcal^{f'} \subset \Pcal^f$ as soon as $f \leq f'$. 

\begin{lemma}\label{lemma:max-front}
The front $f^\star$ defined by 
\[
f^\star_i \coloneqq \max \big\{ j \in [d_i + 1] \colon \forall (x,s) \in \Pcal^\lex \, , \; s_i^1 = \dots = s_i^{j-1} = 0 \big\} \enspace \text{for all} \enspace i \in [m] 
\]
is the greatest front $f$ such that $\Pcal^\lex \subset \Pcal^f$.
\end{lemma}

We call $f^\star$ the \emph{max-front}. \Cref{prop:max-front} below states that the \emph{max-front polyhedron} $\Pcal^{f^\star}$ is a tight approximation of $\Pcal^\lex$, in the sense that $\Pcal^\lex$ stands between $\Pcal^{f^\star}$ and its relative interior. An illustration is given in~\Cref{fig:lex-polyhedron}(b) and (c). We recall that the \emph{affine hull} of a set $S$ is the smallest (inclusionwise) affine subspace containing it, and that the \emph{relative interior} $\relint S$ of $S$ is its interior w.r.t.~the subspace topology on its affine hull. For instance, the relative interior of a point is the point itself. The relative interior of a (closed) line segment between two points is the open line segment between the two points, etc. We note that the relative interior of a nonempty convex set $S$ is always nonempty, and that its closure is equal to the closure of $S$. We refer to~\cite[Chapter~6]{Rockafellar} for a complete account on the topological properties of convex sets.
\begin{proposition}\label{prop:max-front}
We have 
\[
\relint \Pcal^{f^\star} \subset \Pcal^\lex \subset \Pcal^{f^\star} \, .
\] 
Equivalently, $\Pcal^{f^\star}$ is the closure of $\Pcal^\lex$.
\end{proposition}

\begin{figure}[t]
\begin{algorithmic}[1]
\Procedure{ComputeMaxFront}{}
\State $f \coloneqq (1, \dots, 1)$
\Repeat \label{line:loop}
\State $I \coloneqq \varnothing$
\For{$i \in [m]$ such that $f_i \leq d_i$}
\sIf{$\sup \; \{ s_i^{f_i} \colon (x,s) \in \Pcal^f \} \leq 0$}\label{line:LP} $I \coloneqq I \cup \{i\}$
\EndFor
\For{$i \in I$}
\State $f_i \coloneqq f_i + 1$ 
\EndFor
\Until{$I = \varnothing$}\label{line:until}
\State \Return $f$
\EndProcedure
\end{algorithmic}
\caption{Computing the max-front}\label{algo:front}
\end{figure}

We now introduce the algorithm \Call{ComputeMaxFront}{} that computes the front $f^\star$; see \Cref{algo:front}. The principle of the algorithm is the following. It starts from the all-$1$ front, and then repeatedly increment the entries $f_i$ of the current front $f$ as soon as $\sup \; \{ s_i^{f_i} \colon (x,s) \in \Pcal^f \} \leq 0$. The latter test is carried out by solving the linear program 
\[
\text{Maximize} \quad  s_i^{f_i} \quad  \text{subject to} \quad (x, s) \in \Pcal^f \, ,
\]
and checking if the value of this linear program is nonpositive (by convention, this value is $-\infty$ is $\Pcal^f$ if empty). This can be done in polynomial time~\cite[Chapter~13]{Schrijver86}. We obtain the following result:
\begin{theorem}\label{th:max-front}
The algorithm \Call{ComputeMaxFront}{} computes the max-front $f^\star$ in polynomial time, by solving at most $m D$ linear programs.
\end{theorem}

\begin{remark}
As a corollary of \Cref{th:max-front}, we can check in polynomial time if the lex-polyhedron $\Pcal^\lex$ is empty or not. Indeed, by \Cref{prop:max-front}, $\Pcal^\lex$ is nonempty if and only if the max-front polyhedron $\Pcal^{f^\star}$ is nonempty. The latter condition can be verified by first computing $f^\star$, and then testing the emptyness of $\Pcal^{f^\star}$ in polynomial time using linear programming~\cite[Theorem~10.4]{Schrijver86}.
\end{remark}

\begin{example}\label{ex:front}
We consider a lex-polyhedron over the variables $x, y, z$ and defined by the two lex-inequalities $(y,x) \leqlex (x,z)$ and $(x,z) \leqlex (y, y-1)$. It is easy to check that this lex-polyhedron is empty. In standard form, the system of constraints writes as the following set of equalities 
\begin{equation}\label{eq:equalities}
y + s_1^1 = x \, , \enspace x + s_1^2 = z \, , \enspace x + s_2^1 = y \, , \enspace z + s_2^2 = y-1 \, ,
\end{equation}
with the lex-inequalities $(s_1^1, s_1^2) \geqlex 0$ and $(s_2^1, s_2^2) \geqlex 0$. 

We now give the details of the execution of the algorithm \Call{ComputeMaxFront}{}. Initially, the front is given by $f = (1,1)$, and the polyhedron $\Pcal^f$ is defined by the equalities in~\eqref{eq:equalities} and the two inequalities $s_1^1 \geq 0$ and $s_2^1 \geq 0$. The first and third equalities entail that $s_1^1 = -s_2^2$ and so the supremum of $s_1^1$ as well as $s_2^1$ over $\Pcal^f$ is $0$ ($\Pcal^f$ can be checked to be nonempty). Thus, the front is updated to $(2,2)$ after the first iteration of the loop at Line~\lineref{line:loop}, and the updated polyhedron $\Pcal^f$ is now given by the equalities in~\eqref{eq:equalities} together with $s_1^1 = s_2^1 = 0$, and the inequalities $s_1^2 \geq 0$ and $s_2^2 \geq 0$. The equalities $s_1^1 = s_2^1 = 0$ implies that $x = y$, and by the second and fourth equalities of~\eqref{eq:equalities}, we get that $z = x + s_1^2 \geq x = y$, and $y = z + s_2^2 + 1 \geq z + 1$. We deduce that $\Pcal^f$ is empty. In this case, the supremum of $s_1^2$ as well as $s_2^2$ over $\Pcal^f$ is $-\infty$, and so the front is updated to $(3,3)$. This is indeed the max-front, since the original lex-polyhedron is empty.
\end{example}

\subsubsection*{Checking the Strict Feasibility of a Policy.}

We now come back to our original problem, \ie, determining if a policy $\sigma$ is strictly feasible, and computing the associated cell $C_\sigma$ of the congestion diagram. Recall that $\sigma$ is strictly feasible if and only if, for some choice of resources $(r_i^a)$, there is a solution $(\rho, u) \in \R^n \times \R^n$ to the lex-inequality system~\eqref{eq:policy-lex} such that $\rho \geq 0$ and $\rho \neq 0$. This leads to introducing the lex-polyhedron $\Pcal^\lex_\sigma$ over the stationary regime variables $(\rho_i, u_i)$ ($i \in [n]$) and resource variables $(r_i^a)$, defined by the equality and depth-$2$ lex-inequality constraints in the system~\eqref{eq:policy-lex} as well as $n$ depth-$1$ inequalities $\rho_i \geq 0$ for $i \in [n]$. Once turned into standard form, every inequality corresponds to some entry of the max-front $f^\star$. For convenience, we denote by $f^\star_{\sigma_i}$ the entry of corresponding to the inequality $\sigma_i \geq 0$ ($i \in [n]$). By definition of the max-front, the policy $\sigma$ is strictly feasible if and only if $f^\star_{\sigma_i} = 1$ for some $i \in [n]$. Moreover, denoting $\pi$ the projection map onto the resource variables $(r_i^a)$, the cell $C_\sigma$ is precisely the image under $\pi$ of the lex-polyhedron $\Pcal^\lex_\sigma$. We remark $\pi(\Pcal_\sigma^{f^\star})$ is a (closed) convex polyhedron, as a projection of a polyhedron. In consequence, we get the following theorem:
\begin{theorem}\label{th:strict-feasibility}
We can determine in polynomial time whether a policy $\sigma$ is strictly feasible. Moreover, if $\sigma$ is strictly feasible, we have
\[
\relint \pi(\Pcal_\sigma^{f^\star}) \subset C_\sigma \subset \pi(\Pcal_\sigma^{f^\star}) \, ,
\]
or, equivalently, the closure of $C_\sigma$ is given by the convex polyhedron $\pi(\Pcal_\sigma^{f^\star})$.
\end{theorem}

\section{Uniqueness of the Throughput Associated to a Policy}\label{subsec:throughput}

In this section, we consider a fixed policy $\sigma$, and we denote for brevity $P_\tau$ (resp.~$P$ and $\bar P$) the matrix with rows $[P^{\sigma(i)}_\tau]_i$ (resp.~$[P^{\sigma(i)}]_i$ and $[\bar P^{\sigma(i)}]_i$). (We refer to \Cref{subsec:lex-constraints} for the definition of the latter matrices.) Similarly, we denote by $r$ the vector with entries $r_i^{\sigma(i)}$. We consider the affine dynamics determined by the policy $\sigma$, \ie, 
\begin{align}
z_i(t) = r_i + \sum_{j,\tau} (P_\tau)_{ij}z_j(t-\tau) \enspace, \label{e-affine}
\end{align}
and look for stationary solutions $z(t) =u + \rho t$ ($u, \rho \in \R^n$). Equivalently, this amounts to the system 
\begin{equation}\label{eq:hl-policy}
u = r + P u - \bar{P} \, , \qquad \rho = P \rho \enspace, 
\end{equation}
which corresponds to the equality constraints from~\eqref{eq:policy-lex}. We make the following assumption.
\begin{assumption}\label{e-ssimple}
The matrix $P$ has $1$ as an eigenvalue, and this eigenvalue is semisimple.
\end{assumption}
Recall that an eigenvalue is {\em semisimple} if its algebraic and geometric mutiplicities
coincide. 
We next show that under Assumption~\ref{e-ssimple}, the throughput vector~$\rho$ associated to such a stationary solution is unique, and we provide an explicit formula in terms
of the eigenvectors of the matrix $P$.

Denoting by $q$ the multiplicity of the eigenvalue $1$,
and using the semisimplicity assumption, 
we can find a
basis $v^1,\dots,v^q$ of right eigenvectors of $P$, and a basis
$m^1,\dots,m^q$ of left eigenvectors of $P$, so that $Pv^k=v^k$ and
$m^lP=m^l$, and such that $m^kv^l=\delta_{kl}$, for all $k,l\in[q]$,
where $\delta_{kl}$
denotes Kronecker's delta function. See~\cite[\S5,4]{kato} for
background, especially Eq.~(5.36) on p.~41. Then, we define the {\em
  aggregated matrix}
\[
\bar{M}(\tau)  = ( m^k \bar{P}  v^l)_{k,l\in [q]} \in \R^{q\times q}
\]
and the vector $\bar{r}= (m^k r)_{k\in [q]}\in\R^q$.
\begin{theorem}\label{th-unique}
Let $(\rho, u)$ denote any stationary solution of the dynamics~\eqref{e-affine}. 
Then, for generic values of the parameters $\tau\in\timedelays$, the vector $\rho$ is uniquely determined, being given by 
\begin{align}\label{e-np1}
\rho & =\sum_{k\in[q]} \lambda_k v^k \qquad \text{where} \qquad \lambda  = (\lambda_k)_{k\in[q]} = (\bar{M}(\tau))^{-1} \bar{r} \enspace .
\end{align}
\end{theorem}
In~\Cref{sec-mdp}, we show how \Cref{th-unique} extends classical results from semi-Markov decision process theory to the current nonmonotone setting (\ie, where the matrix $P$ may contain negative entries).

\section{Application to Emergency Department Models}\label{sec:experiments}
\subsubsection{Implementation.} The front-advancing algorithm \Call{ComputeMaxFront}{} is implemented in Python. System equations of the form~\eqref{eq:tds} are symbolically defined using SymPy, from which the corresponding lexicographic systems of the form~\eqref{eq:policy-lex} are constructed. Linear programs for checking the front are solved, after assigning numerical values to the time delay and routing proportion parameters, using the HiGHS solver, via the \texttt{linprog} function in the \texttt{scipy.optimize} package. 
Given a strictly feasible policy, we also use SymPy’s \texttt{linsolve} function to symbolically solve the system of linear equations and obtain the associated throughput vector \(\rho\).

The inequalities resulting from the max-front polyhedron are used to generate a congestion phase diagram with the Parma Polyhedra Library (PPL) interface for Python (pyppl). It is then projected onto the resource variables, from which we derive the vertices of the cells.

\subsubsection{The congestion phase diagram revisited.} 
As described in \Cref{sec:model}, we applied the algorithm to the ED model with only junior and senior doctors and cubicles as potentially limiting resources to produce the congestion phase diagram of~\Cref{tab:congestionED}. The matrices induced by the policies of this model satisfy the semisimplicity assumption~\ref{e-ssimple}, implying that the throughputs for each policy is uniquely given as a function of the resources (\Cref{th-unique}), as reported in the third column. The values of the numerical parameters were fixed as follows: \(\pi_{\textrm{care}} = 0.4\), \(\pi_{\mathrm{test1}} = 0.2\), \(\pi_{\mathrm{exit1}} = 0.4\), \(\pi_{\mathrm{test2}} = 0.7\), \({\pi_\mathrm{exit2}} = 0.3\), \(\tau_{\textrm{JC}} = 4\), \(\tau_{\textrm{JS}} = 1\), \(\tau_{\textrm{SC}} = 3\), \(\tau_{\textrm{EC}} = 2\), \(\tau_{\textrm{care}} = 5\), and \(\tau_{\textrm{test}} = 6\).

A few remarks on the congestion phase diagram are in order. Out of the \(32\) total policies, only \(16\) are strictly feasible. Some of these policies correspond to lower dimensional faces of full-dimensional cells, representing special configurations associated with threshold values of the resources. Another remarkable phenomenon is that two distinct policies may yield the same cell and stationary regime $(\rho, u)$. This occurs in the last cell, where it does not make a difference whether senior consultations are bottlenecked by the number of patients or the number of senior doctors, as no senior consultations take place. We interpret this as two bottleneck configurations being equivalent owing to the special structure of the model. 
Finally, two distinct policies may yield the same throughput vector $\rho$, such as in the case of the fifth and seventh cells. Here, despite different bottlenecks for the junior consultations, the lack of senior doctors limiting both for the release of cubicles and the rate of junior-senior synchronization results in the same total throughput. However, we observe in the experiments that these policies have different vectors $u$ (representing ``storage'' levels in each place.) In total, we arrive at a diagram with \(7\) full-dimensional cells, as reported in \Cref{tab:congestionED}.

For the full model of ED, we obtain \(352\) strictly feasible policies out of \(512\) total, yielding \(76\) full-dimensional cells where \(64\) of them are distinct. As an additional test, we have applied our algorithm to the model ``EMS-B'' of emergency call center of Samu de Paris (medical emergency aid) that was presented in~\cite[Figure~5]{boyet2021}. The algorithm returns the same congestion diagram as the one that was previously computed by hand.

\section{Concluding Remarks}

We have developed a general algorithm to compute congestion diagrams of piecewise linear dynamical systems. The heart of this algorithm is a method
to check in polynomial time whether a given policy is strictly feasible
(i.e., corresponds to a non-identically zero throughput), whereas
earlier methods require an exponential time.
A remaining open problem is to develop a polynomial time output sensitive
algorithm to compute all realizable policies, \ie, an algorithm whose execution
time is possibly exponential but that remains polynomially bounded in the number of strictly feasible
policies.
Another open problem is to reinforce~\Cref{th-unique}: this result
shows that every strictly feasible policy determines a unique throughput
vector. However, one may not exclude that for a given allocation of resources,
several policies may be strictly feasible and lead to different values
of $\rho$. This pathology does not appear in the examples
we analyzed; it remains to identify suitable assumptions
implying that it is so.
Finally, another open problem concerns the analysis of nonstationary behaviors:
the conditions for the existence of the limit $\lim_t z(t)/t$ s are still not
understood beyond the monotone case. 

\begin{credits}
  \subsubsection{\ackname}
  
The authors were partially supported by the URGE project of the Bernoulli Lab,
joint between AP-HP \& INRIA. They thank the members of the project, especially Youri Yordanov, Quentin Delannoy, Judith Leblanc, and Christine Fricker.
\end{credits}

\bibliographystyle{splncs04}

\appendix

\section{EMS-B}
For the convenience of the reader we reproduce the timed Petri net with priorities used to model a bi-level emrgency call center with an additional reservoir of mediating responders studied in~\cite{boyet2021}.
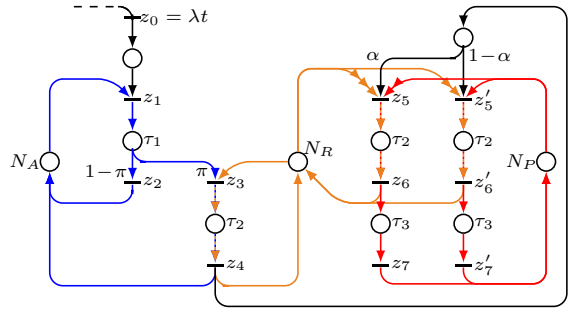
\begin{figure}[h]
  \centering
  \def\tkzscl{.25}
  
  \definecolor{colorARM}{rgb}{0,0,1}
  \definecolor{colorARMres}{rgb}{.92,.5,.11}
  \definecolor{colorAMU}{rgb}{1,0,0}
  \definecolor{colorexit}{rgb}{0.4,0.4,0.4}
  
  \tikzset{place/.style={draw,circle,inner sep=2.5pt,semithick}}
  \tikzset{transition/.style={rectangle, thick,fill=black, minimum width=2mm,inner ysep=0.5pt}}
  \tikzset{token/.style={draw,circle,fill=black!80,inner sep=.35pt}}
  \tikzset{pre/.style={=stealth'}}
  \tikzset{post/.style={->,shorten >=1pt,>=stealth'}}
  \tikzset{-|/.style={to path={-| (\tikztotarget)}}, |-/.style={to path={|- (\tikztotarget)}}}
  \tikzset{Farhi/.style 2 args={dashed,dash pattern=on 1pt off 1pt,#1, postaction={draw,dashed,dash pattern=on 1pt off 1pt,#2,dash phase=1pt}}}
  \tikzset{arrowPetri/.style={>=latex,rounded corners=5pt,semithick}}
  
  \begin{tikzpicture}[scale=\tkzscl,font=\scriptsize]
  
  \def\p{2.2}
  
  \begin{scope}[shift={(0,0)}]
  \node[place] (pool_arm) at ($(-2*\p,-1.5*\p)$) {};
  \node (txt_Na) at ($(pool_arm)+(-1.3,0)$) {$N_A$};
  
  \node[transition]    (q_arrivals)                          at    (0, 0) {};
  \node[place]         (p_inc_calls)                         at    ($(q_arrivals) + (0, \p)$)  {};
  \node[transition]    (q_inc_calls)                         at    ($(p_inc_calls) + (0, \p)$)  {};
  
  \node[place]         (p_arrivals)        at    ($(q_arrivals) + (0,-\p)$) {};
  \node[transition]    (q_debut_NFU)       at    ($(p_arrivals) + (0,-\p)$) {};
  \node[transition]    (q_synchro)         at    ($(p_arrivals) + (2*\p,-\p)$) {};
  
  \node[place]         (pool_res)          at    ($(q_synchro) + (2*\p,.5*\p)$) {};
  \node (txt_Nr) at ($(pool_res)+(1.1,.7)$) {$N_R$};
  
  \draw[->,arrowPetri,colorARM] (p_arrivals)  |- ($(p_arrivals)+(.5*\p,-.5*\p)$) -| (q_synchro);
  \draw[->,arrowPetri,colorARM] (q_debut_NFU) |- ($(q_debut_NFU)+(0,-.5*\p)$) -| (pool_arm);
  \draw[->,arrowPetri,colorARM] (p_arrivals)  -- (q_debut_NFU);
  \draw[->,arrowPetri,colorARM] (q_arrivals) -- (p_arrivals);
  
  \node[place]         (p_synchro)           at    ($(q_synchro)  + (0,-\p)$) {};
  \node[transition]    (q_unsynchro)         at    ($(p_synchro)  + (0,-\p)$) {};
  
  \node[place]      (pool_amu)    at    ($(pool_res)+(6*\p,0)$) {};
  \node (txt_Nm) at ($(pool_amu)+(-1.3,0)$) {$N_P$};
  
  \node[place]         (p_synchro2)          at    ($(pool_amu)  + (-4*\p,.5*\p)$) {};
  \node[transition]    (q_unsynchro2)        at    ($(p_synchro2)  + (0,-\p)$) {};
  \node[transition]    (q_debut_AMU)         at    ($(p_synchro2) + (0,\p)$) {};
  \node[place]         (p_waiting)           at    ($(q_debut_AMU)  + (2*\p,1.5*\p)$) {};
  \node[place]         (p_consult_AMU)       at    ($(q_unsynchro2)+ (0,-1*\p)$) {};
  \node[transition]    (q_end_consult_AMU)   at    ($(p_consult_AMU)+(0,-\p)$) {};
  
  \node[place]         (p_synchro3)          at    ($(p_synchro2)       + (2*\p,0)$) {};
  \node[transition]    (q_unsynchro3)        at    ($(q_unsynchro2)     + (2*\p,0)$) {};
  \node[transition]    (q_debut_AMU_2)       at    ($(q_debut_AMU)      + (2*\p,0)$) {};
  \node[transition]    (q_end_consult_AMU_2) at    ($(q_end_consult_AMU)+ (2*\p,0)$) {};
  \node[place]         (p_consult_AMU_2)     at    ($(p_consult_AMU)    + (2*\p,0)$) {};
  
  \draw[arrowPetri]    (q_inc_calls) |- ($(q_inc_calls)+(-.35*\p,.25*\p)$);
  \draw[dashed, arrowPetri]    (q_inc_calls) |- ($(q_inc_calls)+(-1.5*\p,.25*\p)$);
  \draw[->,arrowPetri] (q_inc_calls) -- (p_inc_calls);
  \draw[->,arrowPetri]                                    (p_inc_calls) -- (q_arrivals);
  \draw[->,arrowPetri,colorARM]                           (q_unsynchro) -- ($(q_unsynchro)+(0,-.5*\p)$)  -|  (pool_arm);
  
  \draw[->,arrowPetri,Farhi={colorARM}{colorARMres}]    (q_synchro) -- (p_synchro);
  \draw[->,arrowPetri,Farhi={colorARM}{colorARMres}]    (p_synchro) -- (q_unsynchro);
  
  \draw[->,arrowPetri,colorARMres]    (pool_res) -- ($(q_synchro)+(.5*\p,.5*\p)$) -- (q_synchro);
  \draw[->,arrowPetri,colorARMres]    (q_unsynchro) |- ($(q_unsynchro)+(.5*\p,-.5*\p)$)-| (pool_res);
  
  \draw[->>>,arrowPetri,colorARMres]   (pool_res) |- ($(q_debut_AMU)+(-.75*\p,.75*\p)$) -- (q_debut_AMU);
  \draw[->,arrowPetri,colorARMres]    (q_unsynchro2) |- ($(q_unsynchro2)+(-\p,-.5*\p)$) -- (pool_res);
  
  \draw[->>,arrowPetri,colorARMres]  (pool_res) |- ($(q_debut_AMU_2)+(-.75*\p,.75*\p)$) -- (q_debut_AMU_2);
  \draw[->,arrowPetri,colorARMres]    (q_unsynchro3) |- ($(q_unsynchro2)+(-\p,-.5*\p)$) -- (pool_res);
  
  \draw[->,arrowPetri,Farhi={colorAMU}{colorARMres}]    (q_debut_AMU) -- (p_synchro2);
  \draw[->,arrowPetri,Farhi={colorAMU}{colorARMres}]    (p_synchro2) -- (q_unsynchro2);
  \draw[->,arrowPetri,Farhi={colorAMU}{colorARMres}]    (q_debut_AMU_2) -- (p_synchro3);
  \draw[->,arrowPetri,Farhi={colorAMU}{colorARMres}]    (p_synchro3) -- (q_unsynchro3);
  
  \draw[->,arrowPetri]    (p_waiting) |- ($(p_waiting)+(-.5*\p,-.5*\p)$) -| (q_debut_AMU);
  \draw[->,arrowPetri]    (p_waiting) -- (q_debut_AMU_2);
  \draw[->,arrowPetri]    (q_unsynchro) |- ($(q_unsynchro)+(1.5*\p,-1*\p)$) -| ($(pool_amu)+(.5*\p,0)$) |- ($(p_waiting)+(0,.75*\p)$) -- (p_waiting);
  
  \draw[->,arrowPetri,colorARM]                           (pool_arm)    |- ($(q_arrivals)+(-.5*\p,.5*\p)$) -- (q_arrivals);
  
  \node (txt_taus) at ($(p_synchro)+(.5*\p,0)$) {$\tau_2$};
  \node (txt_taus) at ($(p_synchro2)+(.5*\p,0)$) {$\tau_2$};
  \node (txt_taus) at ($(p_synchro3)+(.5*\p,0)$) {$\tau_2$};
  \node (txt_tau1) at ($(p_arrivals.center)+(.5*\p,0)$) {${\tau}_1$};
  \node (txt_tau3) at ($(p_consult_AMU.center)+(.5*\p,0)$) {${\tau}_3$};
  \node (txt_tau3) at ($(p_consult_AMU_2.center)+(.5*\p,0)$) {${\tau}_3$};
  \node (txt_pi)  at ($(p_arrivals.center)+(-1.4,-1.65)$) {$1\!-\!\pi$};
  \node (txt_pi3) at ($(p_arrivals.center)+(3.7, -1.65)$){$\pi$};
  \node (txt_alpha3) at ($(p_waiting.center)+(-4.8,-1)$){$\alpha$};
  \node (txt_alpha)  at ($(p_waiting.center) +(1.4,-1)$) {$1\!-\!\alpha$};
  
  \draw[->,arrowPetri,colorAMU]    (q_unsynchro2)      -- (p_consult_AMU);
  \draw[->,arrowPetri,colorAMU]    (p_consult_AMU)     -- (q_end_consult_AMU);
  \draw[->,arrowPetri,colorAMU]    (q_end_consult_AMU) -- ($(q_end_consult_AMU)+(0,-1)$) -| (pool_amu);
  \draw[->>,arrowPetri,colorAMU]   (pool_amu)      |- ($(q_debut_AMU)+(.5*\p,.5*\p)$)      -- (q_debut_AMU);
  
  \draw[->,arrowPetri,colorAMU]    (q_unsynchro3)       -- (p_consult_AMU_2);
  \draw[->,arrowPetri,colorAMU]    (p_consult_AMU_2)     -- (q_end_consult_AMU_2);
  \draw[->,arrowPetri,colorAMU]    (q_end_consult_AMU_2) -- ($(q_end_consult_AMU_2)+(0,-1)$) -| (pool_amu);
  \draw[->,arrowPetri,colorAMU]    (pool_amu)      |- ($(q_debut_AMU_2)+(.5*\p,.5*\p)$)      -- (q_debut_AMU_2);
  
  \node (txt_z0) at ($(q_inc_calls)+(1*\p,0)$) {$z_0 = \lambda t$};
  \node (txt_z1) at ($(q_arrivals)+(.5*\p,0)$) {$z_1$};
  \node (txt_z2) at ($(q_debut_NFU)+(.5*\p,0)$) {$z_2$};
  \node (txt_z3) at ($(q_synchro)+(.5*\p,0)$) {$z_3$};
  \node (txt_z4) at ($(q_unsynchro)+(.5*\p,0)$) {$z_4$};
  \node (txt_z5) at ($(q_debut_AMU)+(.5*\p,0)$) {$z_5$};
  \node (txt_z6) at ($(q_unsynchro2)+(.5*\p,0)$) {$z_6$};
  \node (txt_z7) at ($(q_end_consult_AMU)+(.5*\p,0)$) {$z_7$};
  \node (txt_z5) at ($(q_debut_AMU_2)+(.5*\p,0)$) {$z_5'$};
  \node (txt_z6) at ($(q_unsynchro3)+(.5*\p,0)$) {$z_6'$};
  \node (txt_z7) at ($(q_end_consult_AMU_2)+(.5*\p,0)$) {$z_7'$};
  
  \end{scope}
  
  \end{tikzpicture}	
  
  \caption{Medical emergency call center with a monitored reservoir (EMS-B)~\cite{boyet2021}.}
  \label{fig:SAMU2}
\end{figure}

\section{Sufficient conditions for the existence of stationary regimes}\label{app:suff-cond}
Recall a {\em policy\/} \(\sigma \colon [n]\to \cup_{i\in [n]} \actions_i\) is a map such that \(\sigma(i)\in \actions_i\) for all \(i\).
It is a way of encoding which minimum is chosen in each component \(i\) in~\eqref{eq:tds}. As in~\Cref{subsec:throughput}, to every policy \(\sigma\) we associate a vector \(r^\sigma \in \R^n\), and for each time delay \(\tau\) a matrix \({P_\tau^\sigma \in \R^{n \times n}}\) such that
\begin{equation*}
  r_i^\sigma\coloneqq {(r^{\sigma(i)})}_i,  \qquad {(P_{\tau}^{\sigma})}_{ij} \coloneqq {(P_{\tau}^{\sigma(i)})}_{ij}\,.
\end{equation*}
Specifying a policy \(\sigma\) allows us to define the {\em matrix polynomial}
\[
P^\sigma(\alpha) \coloneqq \sum_\tau P_\tau^\sigma \alpha^\tau \in \R^{n\times n}[\alpha]
\,,
\]
which can be thought of either as a formal polynomial in an {\em indeterminate} $\alpha$, with matrix valued coefficients, or if one prefers, as a matrix-valued polynomial function of the {\em variable} $\alpha$.

We shall make several additional assumptions.
\begin{assumption}\label{as:1}
  For each active policy \(\sigma\) we have that
  \begin{enumerate}[(i)]
    \item For all \(0\leq \alpha <1\), \(1\) is not an eigenvalue of \(P^\sigma(\alpha)\);\label{as:1-1}
    \item The matrix \(P^\sigma(0)\) has no eigenvalue in \([1, +\infty)\).\label{it-1}\label{as:1-2}
  \end{enumerate}
\end{assumption}
The second assumption involves the resolvent of the matrix $P^\sigma(1)$.
We know from~\cite{kato} that if $1$ is a semisimple eigenvalue of the matrix $P^\sigma$,
there exists a sequence
  of matrices \(\LC^\sigma_{-1},\LC^\sigma_0,\LC^\sigma_1,\dots\in \R^{n\times n}\) such that
  \[
  {(I -\alpha P^\sigma(1))}^{-1}  = \frac{\LC_{-1}}{1-\alpha} + \LC_0  + (1-\alpha)\LC_1 + \dots
  \]
  the expansion being absolutely convergent in a sufficiently small
  punctured disk of center \(1\).
  We set
   \[S^\sigma= \sum_{\tau \in \timedelays} (\tau-1) P^\sigma_\tau \enspace .\]
\begin{assumption}\label{as:2}
  For each active policy \(\sigma\),
  \begin{enumerate}[(i)]
    \item Either \(1\) is not an eigenvalue of \(P^\sigma(1)\) or it is a semisimple eigenvalue of \(P^\sigma(1)\);\label{as:2-1}
    \item The matrix \(I + \LC_{-1}^\sigma S^\sigma\) is invertible.
      \label{as:2-2}
  \end{enumerate}
\end{assumption}
It shown in~\cite{HSCC} that under these assumptions,
the system admits a steady state solution of the form:
\begin{equation*}
  x(t) = u + \rho (t + t_1) \quad (u, \rho \in \R^n \, , \; t_1 \geq 0) \,,
  \end{equation*}
  that satisfies
  \begin{equation*}
  u_i + \rho_i (t + t_1) = \min_{a \in \actions_i} \biggl( r_i^a + \sum_{\tau \in \timedelays} \sum_{j \in [n]}(P_{\tau}^a)_{ij} (u_j + \rho_j (t+t_1-\tau)) \biggr) \,,
  \end{equation*}
  for all $t \geq 0$ and $i \in [n]$. Such a solution is known as an \emph{invariant half-line}.
  
\section{Proofs of the Statements of \Cref{subsec:lex-polyhedron,subsec:throughput}}

\begin{proof}[\Cref{lemma:max-front}]
We first prove the inclusion $\Pcal^\lex \subset \Pcal^{f^\star}$. Let $(x,s) \in \Pcal^\lex$. Then $A x + s = b$. Now, let $i \in [m]$. By definition of $f^\star_i$, we have $s^1_i = \dots = s^{f^\star_i-1} = 0$. Since $(s_i^1, \dots, s_i^{d_i}) \geqlex 0$, we must have $s_i^{f^\star_i} \geq 0$ as soon as $f^\star_i \leq d_i$. As a consequence, we have $s^{< f^\star} = 0$ and $s^{f^\star} \geq 0$. This shows that $(x,s) \in \Pcal^{f^\star}$. 

Now, let $f$ be a front such that $\Pcal^\lex \subset \Pcal^f$. Let $i \in [m]$. For all $(x,s) \in \Pcal^\lex$, we have $s_i^1 = \dots = s_i^{f_{i}-1} = 0$ as $(x,s) \in \Pcal^f$. Thus, $f_i \leq f^\star_i$. We deduce that $f \leq f^\star$.\qed
\end{proof}

\begin{proof}[\Cref{prop:max-front}]
The right-hand side inclusion is already proved in~\Cref{lemma:max-front}.

Now, let $i \in [m]$ such that $f^\star_i \leq d_i$. By definition of $f^\star_i$, there exists $(x,s) \in \Pcal^\lex$ such that $s^{f^\star_i}_i \neq 0$. As $\Pcal^\lex \subset \Pcal^{f^\star}$, we have $(x,s) \in \Pcal^{f^\star}$, and so $s^{f^\star_i}_i \geq 0$. Thus, $s^{f^\star_i}_i > 0$. We deduce that no equality $s^{f^\star_i}_i = 0$ ($i \in [m]$, $f^\star_i \leq d_i$) implicitly holds over the set $\Pcal^{f^\star}$.  By~\cite[Section~8.1]{Schrijver86}, this entails that the affine hull of $\Pcal^{f^\star}$ is given by 
\[
\big\{ (x, s) \colon A x + s = b \, , \; s^{< f^\star} = 0 \big\} \, ,
\]
and the relative interior of $\Pcal^{f^\star}$ by 
\[
\relint \Pcal^{f^\star} = \{ (x, s) \colon A x + s = b \, , \; s^{< f^\star} = 0 \, , \, s^{f^\star} > 0 \} \, ,
\]
where $s^{f^\star} > 0$ means that every entry of $s^{f^\star}$ is positive. 
The remaining inclusion $\relint \Pcal^{f^\star} \subset \Pcal^\lex$ straightforwardly follows from the fact that, for all $i \in [m]$, we have $(s_i^1, \dots, s_i^{d_i}) \geqlex 0$, as $s^1_i = \dots = s^{f^\star_i-1} = 0$ and $s^{f^\star_i}_i > 0$ in the case where $f^\star_i \leq d_i$.

The equivalent formulation in terms of the closure is a consequence of~\cite[Corollary~6.3.1]{Rockafellar}, applied to the convex sets $\Pcal^\lex$ and $\Pcal^{f^\star}$, and the fact that $\Pcal^{f^\star}$ is closed.\qed 
\end{proof}

\begin{proof}[\Cref{th:max-front}]
We first deal with the complexity of \Call{ComputeMaxFront}{}. Throughout the execution, the vector $f$ is a front that (strictly) increases at every iteration. As any increasing sequence of fronts has size at most $D+1$ (the $i$th entry of the front can be increased at most $d_i$ times), we deduce that the number of iterations of the loop at Line~\lineref{line:loop} is bounded by $D+1$. There are at most $m$ linear programs solved at every iteration, and $0$ at the last iteration. Thus, at most $m D$ linear programs are solved in total. We deduce that \Call{ComputeMaxFront}{} has polynomial time complexity.

In order to prove the correctness of the algorithm, we make the preliminary observation that, given $i \in [m]$ such that $f_i \leq d_i$, we have $\sup \; \{ s_i^{f_i} \colon (x,s) \in \Pcal^f \} \leq 0$ if and only if $s_i^{f_i} = 0$ for all $(x, s) \in \Pcal^f$. Indeed, let $\gamma \coloneqq \sup \; \{ s_i^{f_i} \colon (x,s) \in \Pcal^f \} \leq 0$. If $s_i^{f_i} = 0$ for all $(x, s) \in \Pcal^f$, then we straightforwardly have $\gamma \leq 0$ (including when $\Pcal_f = \varnothing$, in which case $\gamma = -\infty$). Conversely, suppose that there exists $(x, s) \in \Pcal^f$ such that $s_i^{f_i} \neq 0$. Since $s_i^{f_i} \geq 0$ by definition of $\Pcal^f$, we actually have $s_i^{f_i} > 0$, and so $\gamma > 0$.

Now we prove by induction on the number of iterations of the loop headed at Line~\lineref{line:loop} that $\Pcal^\lex \subset \Pcal^f$. This holds initially when $f$ is the all-$1$ front, since all $(x,s) \in \Pcal^\lex$ must satisfy $A x + s = b$ and $s_i^1 \geq 0$ for each $i \in [m]$. Now suppose that the inclusion $\Pcal^\lex \subset \Pcal^f$ holds at some iteration. We want to prove that $\Pcal^\lex \subset \Pcal^{f'}$, where $f'$ is the front defined by $f'_i \coloneqq f_i+1$ if $i \in I$, $f'_i \coloneqq f_i$ otherwise, and $I$ is the set of $i \in [m]$ such that $\sup \; \{ s_i^{f_i} \colon (x,s) \in \Pcal^f \} \leq 0$. Let $(x,s) \in \Pcal^\lex$. It suffices to prove that $s^{< f'} = 0$ and $s^{f'} \geq 0$, as $(x,s)$ straightforwardly satisfies $A x + s = b$. We distinguish two cases. If $i \notin I$, then $f'_i = f_i$. As  $(x,s) \in \Pcal^\lex \subset \Pcal^f$, we have $s_i^1 = \dots = s_i^{f'_i-1} = 0$, as well as $s_i^{f'_i} \geq 0$ in the case where $f'_i \leq d_i$, by definition of $\Pcal^f$. Now, assume that $i \in I$. As $(x,s) \in \Pcal^\lex \subset \Pcal^f$, we deduce from the observation made in the previous paragraph that $s_i^{f_i} = 0$. Thus, $s_i^1 = \dots = s_i^{f_i - 1} = s^{f_i}_i = 0$. In the case where $f_i < d_i$, the condition $(s_i^1, \dots, s_i^{d_i}) \geqlex 0$ further entails $s_i^{f_i+1} \geq 0$. This shows that we have $s^{< f'} = 0$ and $s^{f'} \geq 0$.

Finally, let $\bar f$ be the front returned by \Call{ComputeMaxFront}{}. Since $\bar f$ satisfies $\Pcal^\lex \subset \Pcal^{\bar f}$ by the previous argument, we deduce that $\bar f \leq f^\star$ thanks to \Cref{lemma:max-front}. It remains to prove $\bar f \geq f^\star$. Let $K \coloneqq \{i \in [m] \colon \bar f_i \leq d_i\} $. If $K = \varnothing$, then for all $i \in [m]$, $\bar f_i = d_i+1$, and so $\bar f \geq f^\star$ trivially holds. Thus, we suppose $K \neq \varnothing$. We observe that for all $i \in K$, there exists a point $z^{(i)} = (x^{(i)},s^{(i)}) \in \Pcal^{\bar f}$ such that $[s^{(i)}]_i^{\bar f_i} > 0$ (otherwise, we would have $\sup \; \{ s_i^{\bar f_i} \colon (x,s) \in \Pcal^{\bar f} \} \leq 0$, which contradicts the stopping condition at Line~\lineref{line:until}). Since $\Pcal^{\bar f}$ is convex, it contains the point $\bar z \coloneqq \frac{1}{|K|} \sum_{i \in K} z^{(i)}$. Thus, decomposing $\bar z = (\bar x, \bar s)$, we have $A \bar x + \bar s = b$, ${\bar s}^{< \bar f} = 0$ and ${\bar s}^{\bar f} \geq 0$. Moreover, for all $k \in K$, 
\begin{align*}
{\bar s}^{\bar f_k}_k & = \frac{1}{|K|} \sum_{i \in K} [s^{(i)}]^{\bar f_k}_k \\
& = \frac{1}{|K|} \Big(\underbrace{[s^{(k)}]^{\bar f_k}_k}_{> 0} + \sum_{i \in K, i \neq k} \underbrace{[s^{(i)}]^{\bar f_k}_k}_{\geq 0} \Big)\\
& > 0 \, .
\end{align*}
To summarize, we have ${\bar s}^{< \bar f} = 0$ and ${\bar s}^{\bar f} > 0$. We deduce that $({\bar s}^1_i, \dots, {\bar s}^{d_i}_i) \geqlex 0$ for all $i \in [m]$, and so $\bar z \in \Pcal^\lex$. Subsequently, for all $i \in K$, the fact that ${\bar s}^{\bar f_i}_i > 0$ implies $\bar f_i \geq f^\star_i$. The inequality $\bar f_i \geq f^\star_i$ still trivially holds when $i \notin K$, since, in this case, $\bar f_i = d_i+1$. We deduce that $\bar f = f^\star$. \qed
\end{proof}

\begin{proof}[\Cref{th:strict-feasibility}]
The first part of the statement is a consequence of \Cref{th:max-front}. The second part is obtained by applying $\pi$ on every side of the relations $\relint \Pcal_\sigma^{f^\star} \subset \Pcal^\lex_\sigma \subset \Pcal_\sigma^{f^\star}$ provided by \Cref{prop:max-front}, and using the fact that $\pi$ and $\relint$ commute for convex sets~\cite[Theorem~6.6]{Rockafellar}. The last equivalent formulation comes from~\cite[Corollary~6.3.1]{Rockafellar}. \qed 
\end{proof}

  \begin{proof}[\Cref{th-unique}]
    We have $\rho = P\rho$. So, $\rho$ can be decomposed
    over the basis of $\ker(I- P)$ given by $v^1,\dots,v^Q$,
so that $\rho=\sum_k \lambda_k v^k$
for some $(\lambda_k)_{k\in[Q]} \in \R^Q$. Substituting
this expression of $\rho$ in $u-Pu=r-\bar{P}\rho$, and left multiplying the latter
equation by $m^k$, for every $k\in [Q]$, to get $\bar{M}(\tau) \lambda = \bar{r}$.

It remains to check that $\bar{M}(\tau)$ is invertible, for generic
values of $\tau$. To see this, consider the special value $\tau\equiv \mathbf{1}$
where $\mathbf{1}$ denotes the matrix whose entries are all equal to $1$.
Then, $(M(\mathbf{1}))_{kl}= m^k P v^l = m^k v^l =\delta_{kl}$, showing
that $M(\mathbf{1})$ is the identity matrix. Then, the determinant
$\det(M(\tau))$, which  is a polynomial in the variables
$\tau_{ij}$, does not vanish at $\tau =\mathbf{1}$. Hence, this determinant
does not vanish for generic values of the variables, and so,
the matrix $M(\tau)$ is generically invertible. \qed 
  \end{proof}
  \section{Comparison with Semi-Markov Decision Processes}\label{sec-mdp}
  When all the coefficients ${(P_\tau^a)}_{ij}$ are nonnegative, and $\sum_{\tau \in\timedelays} \sum_{j\in[n]} (P^a_\tau)_{ij}=1$ for all $i$ and $a \in \actions_i$, the equation~\eqref{eq:tds} coincides with the dynamic programming
  equation of a semi-Markov decision process~\cite{puterman,feinberg1994constrained}. In this context, $z_i(t)$ is interpreted
  as the optimal expected cost starting from state $i$, over an horizon $t$, and
  stationary regimes are known as {\em invariant half-lines}. In the special case
  of Markov decision processes, the delays $\tau$ are identically one. The existence
  of invariant half-lines for models with finite action spaces, and the equivalence
  with the lexicographic formulation~\eqref{eq:stationary-lex} is a classical result in the theory
  of Markov decision processes, see~\cite{dynkin1979controlled,Koh80}.
  The novelty here is to allow the ``probability'' coefficients ${(P_\tau^a)}_{ij}$ to take {\em negative}
  values. 

    \Cref{th-unique} generalizes a classical property of Markov decision processes.
     Indeed, suppose that $P$ is a stochastic matrix and that $\tau=\mathbf{1}$
    (identically one holding times).
    Recall that a $n\times n$ stochastic matrix $P$
    defines a digraph with set of nodes $[n]$ and an arc $i\to j$ if $P_{ij}>0$.
    The {\em final classes} of $P$ are the strongly connected components
    of this digraph that do not have access to any other strongly connected component.
    Let $F_1,\dots,F_q$ denote the set of final classes. The eigenvalue
    $1$ is always semisimple, and its multiplicity coincides with the number
    of final classes. For each final
    class $F_k$, there is a unique invariant measure (left eigenvector
    of $P$ for the eigenvalue $1$, nonnegative and of sum one) $m^k$
    that is supported by $F_k$. Moreover, there is a unique
    eigenvector $v^k$ such that $v^k_{i}=1$ for all $i\in F_k$
    and $v^k_j=0$ for all $j\in \cup_{l\in[q]\setminus \{k\}} F_l$,
    see~\cite[Chap.~8]{bermanandplemmons}, especially Theorems~3.23 and~4.27.
    Then, for $i\in [n]\setminus F_k$, $v^k_i$ gives the probability
    that a trajectory of the Markov chain, starting from $i$, ends
    in the final class $F_k$. Observe that $M(\tau)=I$,
so that $\lambda = \bar{r}$. We shall interpret $r_i$ as an instantaneous payoff received
    in state $i$. Then, $\bar{r}_k = m^k r$ gives the expected
    payoff per time unit, received when the initial state
    is in class $F_k$. Then, \Cref{th-unique} specializes to 
    \[
    \rho_i = \sum_{k\in [q]}v^k_i (m^k r) \enspace,
    \]
    so that $\rho_i$ is the expected mean-payoff per time unit starting from state $i$.
    Formula~\eqref{e-np1} extend this property to the case of ``negative probabilities''.

\end{document}